\documentclass{amsart}

\usepackage{amssymb}

\allowdisplaybreaks

\newtheorem{theo-intro}{Theorem}
\newtheorem{theorem}{Theorem}[section]
\newtheorem{lemma}[theorem]{Lemma}
\newtheorem{proposition}[theorem]{Proposition}
\newtheorem{corollary}[theorem]{Corollary}
\newtheorem{remark}[theorem]{Remark}

\theoremstyle{definition}

\newcommand{\N}{\mathbb{N}}

\newcommand{\C}{\mathbb{C}}

\newcommand{\sig}{\sigma}
\newcommand{\Sig}{\Sigma}
\newcommand{\Degree}{d_{\sigma}}

\newcommand{\EME}{\mathcal{M}}

\newcommand{\Lebesgue}{\mathcal{L}}
\newcommand{\E}{\mathcal{E}}

\newcommand{\rad}{\mathrm{r}}
\newcommand{\gau}{\mathrm{g}}

\newcommand{\dem}{\noindent {\bf Proof. }}
\newcommand{\demone}{\noindent {\bf Proof for $\mathbf{d}_{\Sig}$ bounded. }}
\newcommand{\demtwo}{\noindent {\bf Proof for $\mathbf{d}_{\Sig}$ unbounded. }}
\newcommand{\fin}{\hspace*{\fill} $\square$ \vskip0.2cm}

\newcommand{\8}{\infty}

\newcommand{\kla}{\left ( }
\newcommand{\mer}{\right ) }

\newcommand{\ten}{\otimes}

\newcommand{\pl}{\hspace{.1cm}}

\newcommand{\Si}{\Sigma}

\newcommand{\eps}{\varepsilon}

\newcommand{\M}{{\mathcal M}}

\newcommand{\summ}{\sum\limits}

\begin{document}

\title[Sums of independent noncommutative random
variables] {The norm of sums of independent \\ noncommutative
random variables in $L_p(\ell_1)$}

\author[Junge and Parcet]
{Marius Junge$^{\ast}$ and Javier Parcet$^{\dag}$}

\address{University of Illinois at Urbana-Champaign}

\email{junge@math.uiuc.edu}

\address{Universidad Aut\'{o}noma de Madrid and University of
Illinois at Urbana-Champaign}

\email{javier.parcet@uam.es}

\footnote{$^{\ast}$Partially supported by the NSF DMS-0301116.}
\footnote{$^{\dag}$Partially supported by the Project BFM
2001/0189, Spain.} 

\begin{abstract}
We investigate the norm of sums of independent vector-valued
random variables in noncommutative $L_p$ spaces. This allows us to
obtain a uniform family of complete embeddings of the Schatten
class $S_q^n$ in $S_p(\ell_q^m)$ with optimal order $m \sim n^2$.
Using these embeddings we show the surprising fact that the sharp
type (cotype) index in the sense of operator spaces for $L_p[0,1]$
is $\min(p,p')$ ($\max(p,p')$). Similar techniques are used to
show that the operator space notions of $\mathrm{B}$-convexity and
$\mathrm{K}$-convexity are equivalent.
\end{abstract}

\maketitle

\section*{Introduction}

Sums of independent random variables have a long tradition both in
probability theory and Banach space geometry. More recently, the
noncommutative analogs of these probabilistic results have been
developed \cite{J2,JX,PX} and applied to operator space theory
\cite{J3,J4,PS}. In this paper, we follow this line of research in
studying type and cotype in the sense of operator spaces
\cite{Pa2}. This theory is closely connected to the notions of
$\mathrm{B}$-convexity and $\mathrm{K}$-convexity. Using embedding
results we show that these notions remain equivalent in the
category of operator spaces.

\vskip3pt

We recall from \cite{MauP}  that a Banach space $\mathrm{X}$ is
called $\mathrm{K}$-convex whenever the Gauss projection
$$\mathrm{P}_{\mathbf{G}}: f \in L_2(\Omega; \mathrm{X})
\longmapsto \sum_{k=1}^{\infty} \Big( \int_{\Omega} f(\omega)
\overline{\gau_k(\omega)} \, d\mu(\omega) \Big) \gau_k \in
L_2(\Omega;\mathrm{X})$$ is bounded. Here $\gau_1, \gau_2, \ldots$
are independent standard complex-valued Gaussian random variables
defined over a probability space $(\Omega, \mathsf{A}, \mu)$. An
operator space $\mathrm{X}$ is called $\mathrm{OK}$-convex if
$\mathrm{P}_{\mathbf{G}}$ is completely bounded or equivalently
$S_2(\mathrm{X})$ is $\mathrm{K}$-convex as a Banach space. Using
standard tools from Banach space theory, we know that
$S_2(\mathrm{X})$ is $\mathrm{K}$-convex if and only if
$S_p(\mathrm{X})$ is $\mathrm{K}$-convex for some (any) $1 < p <
\infty$. Therefore this notion does not depend on the parameter
$p$. According to a deep theorem of Pisier \cite{P0}
$\mathrm{K}$-convexity is equivalent to $\mathrm{B}$-convexity.
Following Beck \cite{B}, a Banach space $\mathrm{X}$ is called
$\mathrm{B}$-convex if there exists $n \ge 1$ and $0 < \delta \le
1$ such that $$\frac{1}{n} \inf_{|\alpha_k| = 1} \Big\|
\sum_{k=1}^n \alpha_k x_k \Big\| \le (1 - \delta) \max_{1 \le k
\le n} \|x_k\|$$ holds for any family $x_1, x_2, \ldots, x_n$ of
vectors in $\mathrm{X}$. Giesy proved in \cite{G} that a Banach
space $\mathrm{X}$ is $\mathrm{B}$-convex if and only if
$\mathrm{X}$ does not contain $\ell_1^n$'s uniformly. In the
context of operator spaces, the noncommutative analogue of
$\ell_1^n$ is the Schatten class $S_1^n$, the dual of
$\mathcal{B}(\ell_2^n)$. More generally, one might consider
arbitrary dual spaces $$L_1(\mathcal{A}) = S_1^{n_1} \oplus
S_1^{n_2} \oplus \cdots \oplus S_1^{n_m}$$ of finite dimensional
$\mathrm{C}^*$-algebras. A priori, it is unclear which analogue of
the notion of $\mathrm{B}$-convexity is the right one for operator
spaces. Namely, we could only exclude the $\ell_1^n$'s or all the
$L_1(\mathcal{A})$'s. The following result clarifies this
question.

\begin{theo-intro} \label{Theo-Intro-Pisier}
Let $\mathrm{X}$ be an operator space and let $(\mathcal{A}_n)$ be
a sequence of pairwise different finite dimensional
$\mathrm{C}^*$-algebras. The following are equivalent:
 \begin{enumerate}
  \item[i)] $\mathrm{X}$ is $\mathrm{OK}$-convex.
  \item[ii)] $S_p(\mathrm{X})$ does not contain $\ell_1^n$'s uniformly for
  some {\rm(}any{\rm )} $1<p<\infty$.
  \item[iii)] $S_p(\mathrm{X})$ does not contain $L_1(\mathcal{A}_n)$'s
  uniformly for some {\rm(}any{\rm )} $1<p<\infty$.
 \end{enumerate}
\end{theo-intro}

In fact, we prove a stronger result. We shall say that the spaces
$L_1(\mathcal{A}_n)$'s embed semi-completely uniformly in
$S_p(\mathrm{X})$ when there exists a family of embeddings
$$\Lambda_n: L_1(\mathcal{A}_n) \rightarrow S_p(\mathrm{X})$$
satisfying $\|\Lambda_n\|_{cb} \|\Lambda_n^{-1}\| \le \mathrm{c}$
for some universal constant $\mathrm{c} > 1$. This notion came out
naturally in the paper \cite{Pa}. We refer to \cite{OR} for
further applications of this concept. For the equivalence of ii)
and iii), we prove that if $\ell_1^{n}$ embeds into
$S_p(\mathrm{X})$ with constant $\mathrm{c}_n$, then there is a
map $u: S_1^n \to S_p(\mathrm{X})$ such that $$\big\| u: S_1^n \to
S_p(\mathrm{X}) \big\|_{cb} \big\|u^{-1}: u(S_1^n) \to S_1^n
\big\| \le \mathrm{C} \, \mathrm{c}_n.$$

This map is constructed using (noncommutative) probabilistic
tools. Usually, estimates for sums of noncommutative random
variables are motivated by classical probabilistic inequalities.
Our probabilistic motivation here is given by the following
result. Let us consider a finite collection $f_1, f_2, \ldots,
f_n$ of independent random variables on a probability space
$(\Omega, \mathsf{A}, \mu)$. Then, given $1 \le p < \infty$, the
following equivalence of norms holds

\renewcommand{\theequation}{$\Sigma_p$}
\addtocounter{equation}{-1}

\begin{equation} \label{Equation-S(p,1)}
\Big( \int_{\Omega} \Big[ \sum_{k=1}^n |f_k(\omega)| \Big]^p d
\mu(\omega) \Big)^{1/p} \thicksim \max_{r \in \{1,p\}} \left\{
\Big( \sum_{k=1}^n \int_{\Omega} |f_k(\omega)|^r d \mu(\omega)
\Big)^{1/r} \right\}.
\end{equation}
We shall provide in this paper the natural analog of
(\ref{Equation-S(p,1)}) for noncommutative random variables. This
result requires the use of the so-called \emph{asymmetric $L_p$
spaces}, which will be defined below. Now, going back to the
construction of the map $u: S_1^n \to S_p(\mathrm{X})$, we
consider positive integers $m,n \ge 1$. Then, if $1 \le k \le m$
and $\tau$ stands for the normalized trace, let $\pi_k:
L_p(\tau_n) \rightarrow L_p(\tau_{n^m})$ be the mapping defined by
the relation $\pi_k(x) = 1 \otimes \cdots \otimes 1 \otimes x
\otimes 1 \otimes \cdots \otimes 1$, where $x$ is located at the
$k$-th position and $1$ stands for the identity of $M_n$. Then,
the embedding $u$ can be easily constructed by combining condition
ii) with the following result, which might be of independent
interest.

\renewcommand{\theequation}{\arabic{equation}}

\begin{theo-intro} \label{Theo-Intro-Embedding}
Let $1 < p < \infty$ and $1 \le q \le \infty$. Then, given $n \ge
1$ and $m \ge n^2$, the following map is a complete isomorphism
onto a completely complemented subspace
 \[ x \in S_q^n \longmapsto
 \frac{1}{n^{1/q}} \sum_{k=1}^m \delta_k \otimes \pi_k(x) \in
 L_p(\tau_{n^m};\ell_q^m).\]
Moreover, the $cb$-distance constants are uniformly bounded on
the dimensions.
\end{theo-intro}

Our proof requires $m \ge n^2$ which is different from the
well-known commutative order $m \sim n$. Using type/cotype
estimates we show that the order $m\sim n^2$ is best possible. The
cotype for operator spaces is motivated by the Hausdorff-Young
inequality for non-abelian compact groups. Let $\mathrm{G}$ be a
noncommutative compact group and $\widehat{\mathrm{G}}$ its dual
object. That is, a list of inequivalent irreducible unitary
representations. Given $1 \le p \le 2$, an operator space
$\mathrm{X}$ has Fourier type $p$ with respect to $\mathrm{G}$ if
the $\mathrm{X}$-valued Hausdorff-Young inequality $$\Big(
\sum_{k=1}^n d_k \|A_k\|_{S_{p'}^{d_k}(\mathrm{X})}^{p'}
\Big)^{1/p'} \le_{cb}
\mathcal{K}_p(\mathrm{X},\widehat{\mathrm{G}}) \, \Big(
\int_{\mathrm{G}} \Big\| \sum_{k=1}^n d_k \mbox{tr}(A_k \pi_k(g))
\Big\|_\mathrm{X}^p d\mu(g) \Big)^{1/p}$$ holds for all finite
sequence of matrices $A_1, A_2, \ldots, A_n$ with $A_k \in
M_{d_{\pi_k}} \ten \mathrm{X}$. Here $\mu$ is the normalized Haar
measure and $d_k$ denotes the degree of the irreducible
representation $\pi_k: \mathrm{G} \to U(d_k)$. Moreover, here and
in the following, the symbol $\le_{cb}$ is used to indicate the
corresponding linear map is indeed completely bounded. The notion
of Fourier cotype is dual to the notion of Fourier type stated
above. Following \cite{P2}, we notice that this inequality forces
us to consider an operator space structure on the vector space
where we are taking values. In other words, we need to take values
in operator spaces rather than Banach spaces.

\vskip3pt

Note that the span of the functions of the form $\mbox{tr}(A_k
\pi_k(g))$ is dense in $L_{2}(\mathrm{G})$. In the classical
notion of cotype the right hand side is replaced by a suitable
subset of characters. Since it is not entirely clear which will be
such a canonical subset for arbitrary groups, we follow the
approach of Marcus/Pisier \cite{MP} and consider random Fourier
series of the form $$\sum_{k=1}^n d_k \mbox{tr}(A_k \pi_k(g)
U_{\pi_k})$$ where the $U_{\pi}$'s are random unitaries. However,
the contraction principle allows us to eliminate the coefficients
$\pi(g)$. Therefore, given such a family of random unitaries over
a probability space $(\Omega, \mathsf{A}, \mu)$, a possible notion
of cotype for operator spaces is given by the inequality $$\Big(
\sum_{k=1}^n d_k \|A_k\|_{S_{p'}^{d_k}(\mathrm{X})}^{p'}
\Big)^{1/p'} \le_{cb}
\mathcal{K}_p(\mathrm{X},\widehat{\mathrm{G}}) \, \Big(
\int_{\Omega} \Big\| \sum_{k=1}^n d_k \mbox{tr}(A_k
U_{\pi_k}(\omega)) \Big\|_\mathrm{X}^p d\mu(\omega) \Big)^{1/p}.$$
Although this definition originated from compact groups, in this
formulation only the degrees of the representations of the dual
object and their multiplicity are kept. We may therefore consider
this notion of cotype for arbitrary collections of random
unitaries $(U_{\sigma})$ indexed by $\sig \in \Sigma$ and where
$d_{\sig}$ represents the dimension of $U_{\sig}$. Examples for
$\Si$'s coming from groups are the \emph{commutative} set of
parameters ($\Sig_0 = \N$ and $d_k = 1$ for all $k \ge 1$) which
arises from any non-finite abelian compact group and the set
$\Sig_1 = \N$ with $d_k = k$ for $k \ge 1$, which comes from the
classical Lie group $SU(2)$. As we shall see in this paper, these
two sets of parameters are the most relevant ones in the theory.
Let us mention that random unitaries can also be understood as a
higher dimensional version of random signs or independent
Steinhaus variables. It is rather surprising that in disproving
cotype $q$ larger matrices are not necessarily easier. In part
because the sequence $(d_{\sig})_{\sig \in \Sig}$ provides a new
normalization. Let us also note that Khintchine-Kahane
inequalities are not available in the operator space setting
because they even fail in the level of scalars \cite{Lu,LuP,P2}.
Indeed, type (see Section \ref{Section6} for a definition) and
cotype turn out to be closer notions to Fourier type and cotype
than in the classical theory.

\begin{theo-intro} \label{Theo-Intro-Sharp}
Any infinite dimensional $L_p$ space has:
\begin{itemize}
\item[i)] Sharp $\Sig$-type $\min(p,p')$.
\item[ii)] Sharp $\Sig$-cotype $\max(p,p')$.
\end{itemize}
\end{theo-intro}

The organization of the paper is as follows. Section
\ref{Section1} is devoted to describe the operator space structure
and some basic properties of the asymmetric $L_p$ spaces. These
spaces provide an important tool in this paper. Section
\ref{Section2} contains some preliminary estimates that will be
used in Section \ref{Section3} to prove the analog of
(\ref{Equation-S(p,1)}) for noncommutative random variables. In
Section \ref{Section4} we construct the embedding of $S_q^n$ into
$S_p(\ell_q^m)$ described in Theorem \ref{Theo-Intro-Embedding}.
Section \ref{Section5} is devoted to prove the operator space
version of Pisier's characterization of $\mathrm{K}$-convexity.
Finally, in Section \ref{Section6} we find the sharp operator
space type and cotype indices of $L_p$ spaces.

\section{Asymmetric $L_p$ spaces}
\label{Section1}

Throughout this paper, some basic notions of noncommutative $L_p$
spaces and operator space theory will be assumed, see \cite{P2,P3}
for a systematic treatment. We begin by studying some basic
properties of the asymmetric $L_p$ spaces, defined as follows. Let
$E$ be an operator space and let $\mathcal{M}$ be a semi-finite
von Neumann algebra equipped with a n.s.f. trace $\varphi$. Given
a pair of exponents $2 \le r,s \le \infty$ such that $\frac{1}{p}
= \frac{1}{r} + \frac{1}{s}$, we define the \textbf{asymmetric
$L_p$ space} $L_{(r,s)}(\mathcal{M},\varphi; E)$ as the completion
of $L_p(\mathcal{M}, \varphi) \otimes E$ with respect to the
following norm $$\|x\|_{L_{(r,s)}(\M,\varphi; E)} = \inf_{x =
\alpha y \beta} \, \Big\{ \|\alpha\|_{L_r(\M,\varphi)}
\|y\|_{L_{\infty}(\M,\varphi; E)} \|\beta\|_{L_s(\M,\varphi)}
\Big\},$$ where the infimum runs over all decompositions $x =
\alpha y \beta$ with $\alpha \in L_r(\mathcal{M},\varphi)$, $\beta
\in L_s(\mathcal{M},\varphi)$ and $y \in \mathcal{M}
\otimes_{\mbox{\scriptsize min}} E$. Recall that any
noncommutative $L_p$ space can be realized as $L_p(\M,\varphi;E) =
L_{(2p,2p)}(\M,\varphi;E)$. In this paper, the von Neumann algebra
$\M$ will always be a finite matrix algebra $M_n$ so that the
trace $\varphi$ is unique up to a constant factor. In fact, we
shall only work with the usual trace $\mbox{tr}_n$ of $M_n$ and
its normalization $\tau_n = \frac{1}{n} \mbox{tr}_n$. The spaces
$$S_{(r,s)}^n (E) = L_{(r,s)}(\mbox{tr}_n; E)$$ can be regarded as
the asymmetric version of Pisier's vector-valued Schatten classes.
If $R$ and $C$ stand for the row and column operator Hilbert
spaces, we shall denote by $C_p$ and $R_p$ the interpolation
spaces $[C,R]_{1/p}$ and $[R,C]_{1/p}$. The superscript $n$ will
indicate the $n$-dimensional version. By using elementary
properties of the Haagerup tensor product, it is not difficult to
check that $$S_{(r,s)}^n (E) = C_{r/2}^n \otimes_h E \otimes_h
R_{s/2}^n$$ isometrically. This provides a \emph{natural} operator
space structure for $L_{(r,s)}(\tau_n; E)$.

\begin{lemma} \label{Lemma-CB(Cr,Cs)}
Given $1 \le r,s,t \le \infty$ such that $\frac{1}{s} =
\frac{1}{r} + \frac{1}{t}$, we have $$\begin{array}{l}
\|\alpha\|_{\mathcal{CB} (C_r^n,C_s^n)} = \|\alpha\|_{S_{2t}^n},
\\ \|\beta\|_{\mathcal{CB}(R_r^n,R_s^n)} = \|\beta\|_{S_{2t}^n}.
\end{array}$$
\end{lemma}

\dem Since the row case can be treated similarly, we just prove
the first equality. When $r=s$ the result is trivial. Assume now
that $r = \infty$. We begin by recalling the following well-known
complete isometries $$\mathcal{CB}(C_{\infty}^n,C_s^n) =
R_{\infty}^n \otimes_{\mbox{\scriptsize min}} C_s^n = C_s^n
\otimes_h R_{\infty}^n.$$ Since the Haagerup tensor product
commutes with complex interpolation, we obtain the following
Banach space isometries $$C_s^n \otimes_h R_{\infty}^n =
[C_{\infty}^n \otimes_h R_{\infty}^n, R_{\infty}^n \otimes_h
R_{\infty}^n]_{1/s} = [S_{\infty}^n, S_2^n]_{1/s} = S_{2s}^n.$$
Since $s = t$ when $r = \infty$, the identity holds. Now we take
$s < r < \infty$. In that case, we use the fact that the complex
interpolation space $$S_{2t}^n = [S_{2s}^n,S_{\infty}^n]_{s/r} =
[\mathcal{CB}(C_{\infty}^n,C_s^n),
\mathcal{CB}(C_s^n,C_s^n)]_{s/r} \subset
\mathcal{CB}(C_r^n,C_s^n)$$ is contractively included in
$\mathcal{CB}(C_r^n,C_s^n)$. This gives
$\|\alpha\|_{\mathcal{CB}(C_r^n,C_s^n)} \le
\|\alpha\|_{S_{2t}^n}$. For the lower estimate, we consider the
bilinear form $$\mathcal{CB}(C_r^n,C_s^n) \times
\mathcal{CB}(C_{\infty}^n,C_r^n) \longrightarrow
\mathcal{CB}(C_{\infty}^n,C_s^n),$$ defined by $(\alpha, \beta)
\mapsto \alpha \circ \beta$. Then, recalling that the Banach space
$\mathcal{CB}(C_{\infty}^n,C_p^n)$ is isometrically isomorphic to
the Schatten class $S_{2p}^n$ (see above), we obtain the following
inequality $$\|\alpha \beta\|_{S_{2s}^n} \le
\|\alpha\|_{\mathcal{CB}(C_r^n,C_s^n)} \|\beta\|_{S_{2r}^n}.$$
Taking the supremum over $\beta \in M_n$, we obtain
$\|\alpha\|_{\mathcal{CB}(C_r^n,C_s^n)} \ge
\|\alpha\|_{S_{2t}^n}$. \fin

In the following lemma, we state some basic properties of the
asymmetric $L_p$ spaces which naturally generalize some Pisier's
results in Chapter 1 of \cite{P2}.

\begin{lemma} \label{Lemma-Asymmetric}
The asymmetric Schatten classes satisfy the following properties:
\begin{itemize}
\item[i)] Given $2 \le p,q,r,s,u,v \le \infty$ such
that $\frac{1}{p} = \frac{1}{r} + \frac{1}{u}$ and $\frac{1}{q} =
\frac{1}{s} + \frac{1}{v}$, we have $$\|\alpha x
\beta\|_{S_{(p,q)}^n(E)} \le \|\alpha\|_{S_u^n}
\|x\|_{S_{(r,s)}^n(E)} \|\beta\|_{S_v^n}.$$
\item[ii)] Given $x \in M_n \otimes E$ and $2 \le
p,q \le \infty$, we have $$\|x\|_{M_n(E)} = \sup \Big\{ \|\alpha x
\beta\|_{S_{(p,q)}^n(E)}: \ \|\alpha\|_{S_p^n}, \|\beta\|_{S_q^n}
\le 1 \Big\}.$$ Therefore, any linear map $u: E \rightarrow F$
between operator spaces satisfies $$\|u\|_{cb} = \sup_{n \ge 1} \,
\Big\| id \otimes u: S_{(p,q)}^n(E) \rightarrow S_{(p,q)}^n(F)
\Big\|.$$
\item[iii)] Any block-diagonal matrix
$\mathrm{D}_n(x) \in M_{mn} \otimes E$, with blocks $x_1, x_2,
\ldots, x_n$ in $M_m \otimes E$, satisfies the following identity
$$\|\mathrm{D}_n(x)\|_{S_{(r,s)}^{mn} (E)} = \Big( \sum_{k=1}^n
\|x_k\|_{S_{(r,s)}^m (E)}^p \Big)^{1/p}, \quad \mbox{where} \quad
\textstyle \frac{1}{p} = \frac{1}{r} + \frac{1}{s}.$$
\end{itemize}
\end{lemma}

\dem Let us define $$\alpha \otimes_h id_E \otimes_h
\beta^{\mbox{t}}: \ C_{r/2}^n \otimes_h E \otimes_h R_{s/2}^n
\longrightarrow C_{p/2}^n \otimes_h E \otimes_h R_{q/2}^n$$ to be
the mapping $x \mapsto \alpha x \beta$. Then the inequality stated
in (i) follows by Lemma \ref{Lemma-CB(Cr,Cs)} and the injectivity
of the Haagerup tensor product. Let us prove the first identity in
(ii). We point out that $$\|x\|_{M_n(E)} \ge \sup \Big\{ \|\alpha
x \beta\|_{S_{(p,q)}^n(E)}: \ \|\alpha\|_{S_p^n},
\|\beta\|_{S_q^n} \le 1 \Big\},$$ follows immediately from (i). On
the other hand, following Lemma 1.7 of \cite{P2}, we can write
$\|x\|_{M_n(E)} = \|\alpha_0 x \beta_0\|_{S_1^n(E)}$ for some
$\alpha_0, \beta_0$ in the unit ball of $S_2^n$. Then, we consider
decompositions $\alpha_0 = \alpha_1 \alpha$ and $\beta_0 = \beta
\beta_1$ so that $$\|\alpha_1\|_r = \|\alpha\|_p = 1 = \|\beta\|_q
= \|\beta_1\|_s,$$ with $\frac{1}{p} + \frac{1}{r} = \frac{1}{2} =
\frac{1}{q} + \frac{1}{s}$. Applying (i) one more time, we obtain
$$\|x\|_{M_n(E)} \le \sup \Big\{ \|\alpha x
\beta\|_{S_{(p,q)}^n(E)}: \ \|\alpha\|_{S_p^n}, \|\beta\|_{S_q^n}
\le 1 \Big\}.$$ The second identity in (ii) is immediate. Finally,
we prove (iii). Let $E_k$ be the subspace of $E$ spanned by the
entries of $x_k$. Since $E_k$ is finite dimensional, we can find
$\alpha_k, \beta_k \in M_m$ and $y_k \in M_m \otimes E_k$
satisfying $x_k = \alpha_k y_k \beta_k$ and $$\|x_k\|_{S_{(r,s)}^m
(E)} = \|\alpha_k\|_{S_r^m} \|\beta_k\|_{S_s^m}.$$ By homogeneity,
we may assume that $$\sum_{k=1}^n \|x_k\|_{S_{(r,s)}^m (E)}^p =
\sum_{k=1}^n \|\alpha_k\|_{S_r^m}^r = \sum_{k=1}^n
\|\beta_k\|_{S_s^m}^s.$$ Let us consider the block-diagonal
matrices $\mathrm{D}_n(\alpha)$, $\mathrm{D}_n(y)$ and
$\mathrm{D}_n(\beta)$, made up with blocks $\alpha_k$, $y_k$ and
$\beta_k$ $(1 \le k \le n)$ respectively. The upper estimate
follows by considering the decomposition $\mathrm{D}_n(x) =
\mathrm{D}_n(\alpha) \mathrm{D}_n(y) \mathrm{D}_n(\beta)$. In a
similar way, taking $\rho = 2r/(r-2)$ and $\sigma = 2s/(s-2)$, the
upper estimate also holds for the dual space $S_{(r,s)}^{nm} (E)^*
= S_{(\rho, \sigma)}^{nm} (E^*)$. Therefore, the lower estimate
follows by duality. \fin

\section{Preliminary estimates}
\label{Section2}

In this section we prove the main probabilistic estimates to study
the norm of sums of independent noncommutative random variables.
First a word of notation. Throughout this paper, $e_{ij}$ and
$\delta_k$ will denote the generic elements of the canonical basis
of $M_n$ and $\C^n$ respectively.

\begin{lemma} \label{Lemma-D}
Let $E$ be an operator space and let $\mathrm{D}:
S_1^n(\ell_{\infty}^n(E)) \rightarrow \ell_1^n(E)$ be the mapping
defined by $$\mathrm{D} \Big( \sum_{i,j=1}^n e_{ij} \otimes x_{ij}
\Big) \longmapsto \sum_{k=1}^n \delta_k \otimes x_{kk}^k,$$ where
$x_{ij}^k \in E$ stands for the $k$-th entry of $x_{ij}$. Then
$\mathrm{D}$ is a complete contraction.
\end{lemma}

\dem Let us consider the map $\mathrm{d}: \ell_{\infty}^{n^2} \to
M_n$, defined by $$\mathrm{d} \Big( \sum_{i,j=1}^n \lambda_{ij}
(\delta_i \otimes \delta_j) \Big) = \sum_{k=1}^n \lambda_{kk}
e_{kk}.$$ Since the diagonal projection of $\ell_\8^{n^2}$ onto
$\ell_\8^n$ is a complete contraction, we can use the completely
isometric embedding of $\ell_{\infty}^n$ into the subspace of
diagonal matrices of $M_n$ to deduce that $\mathrm{d}$ is a
complete contraction. If $\mathcal{H}$ stands for a Hilbert space
such that $E$ embeds in $\mathcal{B}(\mathcal{H})$ completely
isometrically, then we consider the mapping $$w:
\mathcal{CB}(\mathcal{B}(\mathcal{H}),\ell_\8^n) \longrightarrow
\mathcal{CB}(\ell_\8^n(\mathcal{B}(\mathcal{H})),
\ell_{\infty}^n(\ell_{\infty}^n))$$ defined by $w(\mathrm{T}) =
id_{\ell_\8^n} \ten \mathrm{T}$. Clearly, $w$ is a complete
contraction. Therefore the linear map $v:
\mathcal{CB}(\mathcal{B}(\mathcal{H}), \ell_\8^n) \rightarrow
\mathcal{CB}(\ell_{\infty}^n(\mathcal{B}(\mathcal{H})), M_n)$,
given by $v(\mathrm{T})= \mathrm{d} \circ w(\mathrm{T})$, is a
complete contraction. Let $S_{\mathcal{H}}^1$ denote the predual
of $\mathcal{B}(\mathcal{H})$. Recalling the completely isometric
embeddings
\begin{eqnarray*}
\ell_\8^n(S_{\mathcal{H}}^1) & \hookrightarrow &
\mathcal{CB}(\mathcal{B}(\mathcal{H}),\ell_\8^n), \\
M_n(\ell_1^n(S_{\mathcal{H}}^1)) & \hookrightarrow &
\mathcal{CB}(\ell_\8^n(\mathcal{B}(\mathcal{H})), M_n),
\end{eqnarray*}
we deduce that $u: \ell_\8^n(S_{\mathcal{H}}^1) \to
M_n(\ell_1^n(S_{\mathcal{H}}^1))$, defined by the relation $$u
\Big( \sum_{k=1}^n \delta_k \otimes a_k \Big) = \sum_{k=1}^n
e_{kk} \otimes (\delta_k \otimes a_k),$$ is a complete
contraction. Namely, given $a \in
\ell_{\infty}^n(S_{\mathcal{H}}^1)$ let $$\mathrm{T}_a(x) =
\sum_{k=1}^n \mbox{tr}(a_k^{\mathrm{t}}x) \delta_k \in
\mathcal{CB}(\mathcal{B}(\mathcal{H}), \ell_{\infty}^n).$$ Then,
it can be easily checked that $u(a) = v(\mathrm{T}_a)$. Hence, $u$
can be regarded as the restriction of $v$ to
$\ell_{\infty}^n(S_{\mathcal{H} }^1)$. This proves that $u$ is a
complete contraction. On the other hand, the original map
$\mathrm{D}$ is now given by the restriction of the adjoint map
$u^*: S_1^n(\ell_{\infty}^n(\mathcal{B}(\mathcal{H}))) \rightarrow
\ell_1^n(\mathcal{B}(\mathcal{H}))$ to the subspace
$S_1^n(\ell_{\infty}^n(E))$. \fin

In the following result we shall need the following description of
the Haagerup tensor norm. Given an operator space $E$, we denote
by $M_{p,q}(E)$ the space of $p \times q$ matrices with entries in
$E$. The norm in $M_{p,q}(E)$ is given by embedding it into the
upper left corner of $S_{\infty}^n(E)$ with $n = \max(p,q)$. Now,
for any pair $E_1, E_2$ of operator spaces, let $x_1 \in
M_{p,m}(E_1)$ and $x_2 \in M_{m,q}(E_2)$. We will denote by $x_1
\odot x_2$ the matrix $x$ in $M_{p,q}(E_1 \otimes E_2)$ defined by
$$x(i,j) = \sum_{k=1}^m x_1(i,k) \otimes x_2(k,j).$$ Then, given a
family $E_1, E_2, \ldots, E_n$ of operator spaces and given $$x
\in M_m \otimes \big( E_1 \otimes E_2 \otimes \cdots \otimes E_n
\big),$$ we define the norm of $x$ in the space $S_{\infty}^m
\big( E_1 \otimes_h E_2 \otimes_h \cdots \otimes_h E_n \big)$ as
follows $$\|x\|_m = \inf \left\{ \prod_{k=1}^{n+1}
\|x_k\|_{M_{p_k,p_{k+1}}(E_k)} \, \Big| \ p_1 = m = p_{n+1}
\right\},$$ where the infimum runs over all possible
decompositions $$x = x_1 \odot x_2 \odot \cdots \odot x_n \qquad
\mbox{with} \qquad x_k \in M_{p_k,p_{k+1}}(E_k).$$

\begin{lemma} \label{Lemma-Matrix}
Let $E$ be an operator space and $x_1, x_2, \ldots, x_n \in M_m
\otimes E$. If there are elements $a_k(r), b_k(s) \in M_m$ and
$y_k(r,s) \in M_m \otimes E$ with $1 \le r \le \rho$ and $1 \le s
\le \sigma$ such that $$x_k = \sum_{r,s} a_k(r) y_k(r,s) b_k(s)$$
holds for $1 \le k \le n$, then we have
\begin{eqnarray*}
\lefteqn{\Big\| \sum_{k=1}^n \delta_k \otimes x_k
\Big\|_{L_p(\tau_m; \ell_1^n(E))}} \\ & \le & \Big\| \sum_{k,r}
a_k(r)a_k(r)^* \Big\|_p^{1/2} \sup_k \Big\| \Big( \, y_k(r,s) \,
\Big) \Big\|_{\infty} \Big\| \sum_{k,s} b_k(s)^*
b_k(s)\Big\|_p^{1/2}.
\end{eqnarray*}
\end{lemma}

\dem Let us consider the positive matrices $a,b \in M_m$ defined
by $$a = \Big( \sum_{k,r} a_k(r) a_k(r)^* \Big)^{1/2} \quad
\mbox{and} \quad b = \Big( \sum_{k,s} b_k(s)^* b_k(s)
\Big)^{1/2}.$$ Then, we can find matrices $\alpha_k(r)$ and
$\beta_k(s)$ in $M_m$ satisfying
\begin{eqnarray*}
a_k(r) & = & a \alpha_k(r), \\ b_k(s) & = & \beta_k(s) b,
\end{eqnarray*}
and such that
\begin{eqnarray*}
\Big\| \sum_{k,r} \alpha_k(r) \alpha_k(r)^* \Big\|_{\infty} & \le
& 1, \\ \Big\| \sum_{k,s} \beta_k(s)^* \beta_k(s) \Big\|_{\infty}
& \le & 1.
\end{eqnarray*}
Let us define, for $1 \le k \le n$, the matrices $$z_k =
\sum_{r,s} \alpha_k(r) y_k(r,s) \beta_k(s).$$ Then, by the
definition of $L_p(\tau_m; \ell_1^n(E))$, we have $$\Big\|
\sum_{k=1}^n \delta_k \otimes x_k \Big\|_p \le \Big\| \sum_{k,r}
a_k(r)a_k(r)^* \Big\|_p^{1/2} \Big\| \sum_{k=1}^n \delta_k \otimes
z_k \Big\|_{\infty} \Big\| \sum_{k,s} b_k(s)^*
b_k(s)\Big\|_p^{1/2},$$ where $\| \ \|_{\infty}$ denotes here the
norm on $´S_{\infty}^m(\ell_1^n(E))$. Therefore, it suffices to
estimate the middle term on the right. To that aim, we consider
the matrices
\begin{eqnarray*}
\alpha & = & (\cdots, \alpha_k(r) \otimes e_{1k}, \cdots) \ \in
M_{m,mn \rho}(R_{\infty}^n), \\ \beta & = & (\cdots, \beta_k(s)
\otimes e_{k1}, \cdots)^{\mbox{t}} \in M_{mn
\sigma,m}(C_{\infty}^n).
\end{eqnarray*}
Moreover, if $y_k = \big( y_k(r,s) \big) \in M_{m \rho,m \sigma}
\otimes E$, we also consider the matrix $$y = \sum_{k=1}^n e_{kk}
\otimes (\delta_k \otimes y_k) \in M_{mn \rho,mn
\sigma}(\ell_{\infty}^n(E)).$$ Finally, we notice that
$$\sum_{k=1}^n \delta_k \otimes z_k = (id \otimes
\mathrm{D})(\alpha \odot y \odot \beta).$$ In particular, Lemma
\ref{Lemma-D} gives $$\Big\| \sum_{k=1}^n \delta_k \otimes z_k
\Big\|_{\infty} \le \Big\| \sum_{k,r} \alpha_k(r) \alpha_k(r)^*
\Big\|_{\infty}^{1/2} \sup_k \Big\| \Big( \, y_k(r,s) \, \Big)
\Big\|_{\infty} \Big\| \sum_{k,s} \beta_k(s)^* \beta_k(s)
\Big\|_{\infty}^{1/2}. $$ This yields the assertion since the
first and third terms on the right are $\le 1$. \fin

Let us consider two positive integers $l$ and $n$. Then, given $1
\le p \le \infty$, we define $\pi_k: L_p(\tau_l) \rightarrow
L_p(\tau_{l^n})$ for each $1 \le k \le n$ to be the mapping
defined by the relation $$\pi_k(x) = 1 \otimes \cdots \otimes 1
\otimes x \otimes 1 \otimes \cdots \otimes 1,$$ where $x$ is
located at the $k$-th position and $1$ stands for the identity of
$M_l$. These operators will appear quite frequently throughout
this paper. In the following lemma, we give an upper estimate of
the $L_p$ norm of certain sums of positive matrices constructed by
means of the mappings $\pi_k$.

\begin{lemma} \label{Lemma-Rosenthal}
Let $a^1, a^2, \ldots, a^n \in M_{ml}$ be a collection of positive
matrices and let $1 \le p < \infty$. Then, we have $$\Big\|
\sum_{k=1}^n \pi_k (a^k) \Big\|_{L_p(\tau_{ml^n})} \le c p \, \max
\left\{ \Big\| \frac{1}{l} \sum_{k=1}^n \sum_{i=1}^l a_{ii}^k
\Big\|_{L_p(\tau_m)}, \Big( \sum_{k=1}^n
\|a^k\|_{L_p(\tau_{ml})}^p \Big)^{1/p} \right\}.$$ Here $c$
denotes an absolute constant not depending on the dimensions.
\end{lemma}

\dem By homogeneity, we may assume that the maximum on the right
is $1$. Let $\E$ be the conditional expectation onto
$L_p(\tau_m)$, regarded as a subspace of $L_p(\tau_{ml^n})$. Let
$b_k$ stand for $\pi_k (a^k)$, by the triangle inequality $$\Big\|
\sum_{k=1}^n \pi_k (a^k) \Big\|_p \le \Big\| \sum_{k=1}^n \E (b_k)
\Big\|_{L_p(\tau_m)} + \ \Big\| \sum_{k=1}^n b_k - \E(b_k)
\Big\|_p = \mathrm{A} + \mathrm{B}.$$ Let us observe that
$$\E(b_k) = \tau_l(a^k) = \frac{1}{l} \sum_{i=1}^l a_{ii}^k.$$
Therefore, by assumption, the term A may be estimated by $1$. To
estimate B, we first assume that $p>2$. Then, applying the
Burkholder inequality given in \cite{JX} for noncommutative
martingales, we obtain
\begin{eqnarray*}
\mathrm{B} & \le & c p \, \max\left\{ \Big\| \sum_{k=1}^n
\E(b_k^2)-\E(b_k)^2 \Big\|_{L_{p/2}(\tau_m)}^{1/2}, \Big(
\sum_{k=1}^n \big\| b_k - \E (b_k) \big\|_p^p \Big)^{1/p} \right\}
\\ & \le & c p \, \max\left\{ \Big\| \sum_{k=1}^n \E(b_k^2)
\Big\|_{L_{p/2}(\tau_m)}^{1/2}, \ 2 \ \Big( \sum_{k=1}^n
\|a^k\|_{L_p(\tau_{ml})}^p \Big)^{1/p} \right\}.
\end{eqnarray*}
On the other hand, since $q = p/2 > 1$, we invoke Lemma 5.2 of
\cite{JX} to obtain
\begin{eqnarray*}
\Big\| \sum_{k=1}^n \E(b_k^2) \Big\|_{L_{p/2}(\tau_m)} & \le &
\Big\| \sum_{k=1}^n \E(b_k)
\Big\|_{L_{2q}(\tau_m)}^{\frac{2(q-1)}{2q-1}} \Big( \sum_{k=1}^n
\|b_k\|_{2q}^{2q} \Big)^{\frac{1}{2q-1}} \\ & = & \Big\|
\sum_{k=1}^n \E(b_k) \Big\|_{L_{2q}
(\tau_m)}^{\frac{2(q-1)}{2q-1}} \Big( \sum_{k=1}^n \|
a^k\|_{L_p(\tau_{ml})}^p \Big)^{\frac{1}{2q-1}}
\end{eqnarray*}
The first factor on the right is a power of A. By assumption, the
second factor may be estimated by 1. When $1 \le p \le 2$, we
proceed in a different way. Given a subset $\Gamma$ of $\{1,2,
\ldots, n\}$ with cardinal $|\Gamma|$, let us consider the
conditional expectation $\E_{\Gamma}$ onto
$L_p(\tau_{ml^{|\Gamma|}})$ given by $$\E_{\Gamma} \Big(z_0
\otimes \bigotimes_{k=1}^n z_k \Big) = \prod_{k \notin \Gamma}
\tau_l(z_k) \Big( z_0 \otimes \bigotimes_{k \in \Gamma} z_k
\Big),$$ with $z_0 \in M_m$ and $z_k \in M_l$ for $1 \le k \le n$.
Therefore, since $x_k = b_k - \E (b_k)$ are independent mean $0$
random variables, we can estimate B as follows. For any family of
signs $\varepsilon_k = \pm 1$ with $1 \le k \le n$, let $\Gamma =
\{k: \, \varepsilon_k = 1\}$. Then we have $$\Big\| \sum_{k=1}^n
x_k \Big\|_p \le \Big\| \E_{\Gamma} \Big( \sum_{k =1}^n
\varepsilon_k x_k \Big) \Big\|_p + \Big\| \E_{\Gamma^{\mbox{\tiny
c}}} \Big( \sum_{k =1}^n \varepsilon_k x_k \Big) \Big\|_p \le 2
\Big\| \sum_{k =1}^n \varepsilon_k x_k \Big\|_p.$$ Then, if we
write $\mathrm{r}_1, \mathrm{r}_2, \ldots$ for the classical
Rademacher variables, we use the fact that $L_p(\tau_{ml^n})$ has
Rademacher type $p$ to obtain
\begin{eqnarray*}
\mathrm{B} = \Big\| \sum_{k=1}^n x_k \Big\|_p & \le & 2 \Big(
\int_0^1 \Big\| \sum_{k=1}^n \mathrm{r}_k(t) x_k \Big\|_p^2 dt
\Big)^{1/2} \\ & \le & 2 \Big( \sum_{k=1}^n \|x_k\|_p^p
\Big)^{1/p}
\\ & \le & 4 \Big( \sum_{k=1}^n \|a^k\|_{L_p(\tau_{ml})}^p
\Big)^{1/p}.
\end{eqnarray*}
This yields the assertion for $1 \le p \le 2$. Therefore, the
proof is completed. \fin

\section{Proof of noncommutative (\ref{Equation-S(p,1)})}
\label{Section3}

In this section, we prove the noncommutative analog of the
equivalence of norms (\ref{Equation-S(p,1)}) described in the
Introduction. However, before that we need to set some notation.
Let $E$ and $F$ be operator spaces such that $(E,F)$ is a
compatible pair for interpolation. In what follows, we shall
denote by $J_t(E,F)$ and $K_t(E,F)$ the $J$ and $K$ functionals on
$(E,F)$ endowed with their natural operator space structures as
defined in \cite{X}. Moreover, given $2 \le r,s \le \infty$, we
shall write $\ell_{(r,s)}^n(E)$ to denote the linear space $E^n$
endowed with the operator space structure which arises from the
natural identification with the diagonal matrices of $S_{(r,s)}^n
(E)$. Then, given $1 \le p,q \le \infty$, we use these spaces to
define the operator space $$J_{p,q}^n (M_l;E) = \bigcap_{r,s \in
\{2p,2q\}} \ell_{(r,s)}^n \big( L_{(r,s)} (\tau_l; E) \big).$$

\begin{lemma} \label{Lemma-Projection}
Let $1 \le p,q \le \infty$ and $\lambda_p = l^{-1/p}$. Then, given
$t = l^{\frac{1}{2p} - \frac{1}{2q}}$, we consider the mapping $u:
J_{p,q}^n(M_l;E) \rightarrow J_t(C_p^{nl},C_q^{nl}) \otimes_h E
\otimes_h J_t(R_p^{nl},R_q^{nl})$, defined by the relation $$u
\Big( \sum_{k=1}^n \delta_k \otimes x_k \Big) = \lambda_p
\sum_{k=1}^n e_{kk} \otimes x_k.$$ Then, $u$ is a complete
isometry with completely contractively complemented image.
\end{lemma}

\dem By the injectivity of the Haagerup tensor product, it can be
checked that $$J_t(C_p^{nl},C_q^{nl}) \otimes_h E \otimes_h
J_t(R_p^{nl},R_{q}^{nl}) = \bigcap_{r,s \in \{2p,2q\}}
\lambda_p^{-1} L_{(r,s)} (\mbox{tr}_n \otimes \tau_l; E).$$ Taking
diagonals at both sides, the first assertion follows. In order to
see that the diagonal is completely contractively complemented, we
use the standard diagonal projection $$\mathrm{P}((x_{ij})) =
\int_{\{-1,1\}^n} (\eps_ix_{ij}\eps_j) \, d \mu(\eps),$$ where
$\mu$ is the normalized counting measure on $\{-1,1\}^n$. Here
$x=(x_{ij})$ is a $n\times n$ matrix with entries in $M_l \otimes
E$. By means of the second identity of Lemma
\ref{Lemma-Asymmetric} (ii), it clearly suffices to check that
$\mathrm{P}$ is a contraction on $L_{(r,s)} (\mbox{tr}_n \otimes
\tau_l; E)$. For each $\varepsilon \in \{-1,1\}^n$, we consider
the matrix $a_{\varepsilon} = (\varepsilon_i \delta_{ij})$. Then
we have $$\|\mathrm{P}(x)\|_{(r,s)} = \Big\| 2^{-n}
\sum_{\varepsilon} a_{\varepsilon} x a_{\varepsilon}
\Big\|_{(r,s)} \le 2^{-n} \sum_{\varepsilon} \|a_{\varepsilon} x
a_{\varepsilon}\|_{(r,s)} = \|x\|_{(r,s)},$$ since
$a_{\varepsilon}$ is unitary for any $\varepsilon \in \{-1,1\}^n$.
This completes the proof. \fin

\begin{proposition} \label{Proposition-Lambda-Cap}
Let $E$ be an operator space and let $1 \le p < \infty$. Then
$$\Lambda_{p1}: \sum_{k=1}^n \delta_k \otimes x_k \in J_{p,1}^n
(M_l;E) \longmapsto \sum_{k=1}^n \delta_k \otimes \pi_k(x_k) \in
L_p(\tau_{l^n}; \ell_1^n(E))$$ is a completely bounded map with
$\|\Lambda_{p1}\|_{cb} \le cp$, where $c$ is independent of $l$
and $n$.
\end{proposition}

\dem Given $t = l^{\frac{1}{2p} - \frac{1}{2}}$, if we regard
$J_t(C_p^{nl},R_{\infty}^{nl}) \otimes_h E \otimes_h
J_t(R_p^{nl},C_{\infty}^{nl})$ as a space of $n \times n$ matrices
with entries in $M_l \otimes E$, then we consider its diagonal
subspace $\mathcal{J}_{p,1}^n(M_l;E)$. By Lemma
\ref{Lemma-Projection}, it suffices to check that the mapping
$$\widetilde{\Lambda}_{p1}: \sum_{k=1}^n e_{kk} \otimes x_k \in
\mathcal{J}_{p,1}^n(M_l;E) \longmapsto \sum_{k=1}^n \delta_k
\otimes \pi_k(x_k) \in L_p(\tau_{l^n}; \ell_1^n(E)),$$ satisfies
$\|\widetilde{\Lambda}_{p1}\|_{cb} \le cp \, \lambda_p$. Given $m
\ge 1$, let us consider $x \in S_p^m(\mathcal{J}_{p,1}^n(M_l;E))$
of norm less than one. Let us consider the spaces
\begin{eqnarray*} F_1 & = & C_p^m \otimes_h J_t(C_p^{nl},R_{\infty}^{nl})
\\ F_2 & = & J_t(R_p^{nl},C_{\infty}^{nl}) \otimes_h R_p^m.
\end{eqnarray*}
Since the space $S_p^m(\mathcal{J}_{p,1}^n(M_l;E))$ embeds
completely isometrically in $F_1 \otimes_h E \otimes_h F_2$, we
can write $x = a \odot y \odot b$, with $a \in M_{1,mln}(F_1)$, $y
\in M_{mln}(E)$, $b \in M_{mln,1}(F_2)$ so that
$\|y\|_{M_{mnl}(E)} < 1$ and
\begin{eqnarray*}
\max \Big\{ \|a\|_{S_{2p}^{mln}}, \, t \|a\|_{C_p^m \ten_h
R_{\infty}^{ml^2n^2}} \Big\} & < & 1 \\ \max \Big\{ \hskip1pt
\|b\|_{S_{2p}^{mln}} \hskip1pt, \, t \|b\|_{C_{\infty}^{ml^2n^2}
\ten_h R_p^m} \Big\} & < & 1.
\end{eqnarray*}
Now, if we write $a = (a_{ij})$, $y = (y_{ij})$ and $b = (b_{ij})$
as $n \times n$ matrices of $ml \times ml$ matrices, we have $$x_k
= \sum_{i,j=1}^n a_{ki} y_{ij} b_{jk} \quad \mbox{where} \quad x =
\sum_{k=1}^n e_{kk} \otimes x_k.$$ Therefore, we have
$$\widetilde{\Lambda}_{p1}(x) = \sum_{k=1}^n \sum_{i,j=1}^n
\delta_k \otimes \pi_k(a_{ki}) \pi_k(y_{ij}) \pi_k(b_{jk}).$$
According to Lemma \ref{Lemma-Matrix}, we deduce
\begin{eqnarray*}
\big\| \widetilde{\Lambda}_{p1}(x) \big\|_p & \le & \Big\|
\sum_{k,i} \pi_k (a_{ki} a_{ki}^*) \Big\|_p^{1/2} \sup_k \Big\|
\Big( \, \pi_k(y_{ij}) \, \Big) \Big\|_{\infty} \Big\| \sum_{k,j}
\pi_k(b_{jk}^* b_{jk}) \Big\|_p^{1/2},
\end{eqnarray*}
where $\|\widetilde{\Lambda}_{p1}(x)\|_p$ stands for the norm of
$\widetilde{\Lambda}_{p1}(x)$ in
$S_p^m(L_p(\tau_{l^n};\ell_1^n(E)))$. As we know, the middle term
on the right is bounded above by $1$. On the other hand, Lemma
\ref{Lemma-Rosenthal} allows us to write
\begin{eqnarray*}
\lefteqn{\Big\| \sum_{k,i} \pi_k (a_{ki} a_{ki}^*)
\Big\|_{L_p(\mathrm{tr}_m \otimes \tau_{l^n})}} \\ & \le & cp \,
\max \left\{ \Big\| \sum_{k,i} \E(a_{ki} a_{ki}^*) \Big\|_{S_p^m},
\, \Big( \sum_{k=1}^n \Big\| \sum_{i=1}^n a_{ki} a_{ki}^*
\Big\|_{L_p(\mathrm{tr}_m \otimes \tau_l)}^p \Big)^{1/p} \right\}
\\ & \le & cp \, \max \left\{ \frac{1}{l} \Big\| \sum_{r,s=1}^l
\sum_{k,i} a_{ki}(r,s) a_{ki}(r,s)^* \Big\|_{S_p^m}, \, \lambda_p
\|a a^*\|_{S_p^{nml}} \right\}.
\end{eqnarray*}
The last inequality follows from the fact that the projection onto
block diagonal matrices is completely contractive, see Corollary
1.3 of \cite{P2}. Now, recalling that $$\Big\| \sum_{r,s=1}^l
\sum_{k,i} a_{ki}(r,s) a_{ki}(r,s)^* \Big\|_{S_p^m} = \|a\|_{C_p^m
\otimes_h R_{\infty}^{ml^2n^2}}^2,$$ we obtain $$\Big\| \sum_{k,i}
\pi_k (a_{ki} a_{ki}^*) \Big\|_{L_p(\mathrm{tr}_m \otimes
\tau_{l^n})} \le cp \lambda_p.$$ Since the same estimate holds for
$b$, we get the desired estimate for
$\|\widetilde{\Lambda}_{p1}\|_{cb}$. \fin

Proposition \ref{Proposition-Lambda-Cap} provides an upper
estimate for the norm of sums of independent noncommutative random
variables in $L_p(\ell_1^n(E))$. Now, we are interested on the
dual version of this result. Hence, it is natural to consider the
operator spaces $K_{p,q}^n(M_l;E)$, which arise when replacing
intersections by sums in $J_{p,q}^n(M_l;E)$. That is, we define
$$K_{p,q}^n (M_l;E) = \sum_{r,s \in \{2p,2q\}} \ell_{(r,s)}^n
\big( L_{(r,s)} (\tau_l; E) \big).$$

\begin{remark} \label{Remark-K-Diagonal}
\emph{Arguing as in Lemma \ref{Lemma-Projection}, we can regard
$K_{p,q}^n(M_l;E)$ as the diagonal in the Haagerup tensor product
$K_t(C_p^{nl},C_q^{nl}) \otimes_h E \otimes_h K_t(R_p^{nl},
R_q^{nl})$ normalized by $\lambda_p = l^{-1/p}$. The projection
$\mathrm{P}$ onto the diagonal is also a complete contraction.}
\end{remark}

\begin{lemma} \label{Lemma-Contractions}
Let $1 \le p \le q \le \infty$ and let $E$ be an operator space.
Then, given a positive integer $n \ge 1$, the following identity
maps are complete contractions $$\begin{array}{l} id:
\ell_{(2p,2q)}^n (E) \longrightarrow \ell_q^n (E) \\ id:
\ell_{(2q,2p)}^n (E) \longrightarrow \ell_q^n (E). \end{array}$$
\end{lemma}

\dem Given $2 \le r,s \le \infty$, we can argue as in Lemma
\ref{Lemma-Projection} to see that the diagonal projection
$\mathrm{P}: S_{(r,s)}^n(E) \rightarrow \ell_{(r,s)}^n (E)$ is a
complete contraction. Therefore, the complex interpolation space
between the diagonals of two asymmetric Schatten classes is the
diagonal of the interpolated asymmetric Schatten class. That is,
we have the following complete isometries $$\begin{array}{rcl}
\ell_{(2p,2q)}^n(E) = \big[ \ell_q^n(E), \ell_{(2,2q)}^n(E)
\big]_{\theta} & \quad & \ell_{(2,2q)}^n(E) = \big[
\ell_{(2,\infty)}^n(E), \ell_1^n(E) \big]_{\gamma} \\
\ell_{(2q,2p)}^n(E) = \big[ \ell_q^n(E), \ell_{(2q,2)}^n(E)
\big]_{\theta} & \quad & \ell_{(2q,2)}^n(E) = \big[
\ell_{(\infty,2)}^n(E), \ell_1^n(E) \big]_{\gamma}
\end{array}$$
completely isometrically for $\frac{1}{2p} = \frac{1-\theta}{2q} +
\frac{\theta}{2}$ and $\gamma = 1/q$. Hence, it suffices to show
the result for $p=1$ and $q = \infty$. That is, we have to see
that the identity mappings $$\begin{array}{l} id:
\ell_{(2,\infty)}^n (E) \longrightarrow \ell_{\infty}^n (E) \\ id:
\ell_{(\infty,2)}^n (E) \longrightarrow \ell_{\infty}^n (E),
\end{array}$$ are complete contractions. In other words, we have
to consider the diagonals of $R_{\infty}^n \otimes_h E \otimes_h
R_{\infty}^n$ and $C_{\infty}^n \otimes_h E \otimes_h
C_{\infty}^n$. However, we recall the completely isometric
isomorphisms $E \otimes_h R_{\infty}^n = E
\otimes_{\mbox{\scriptsize min}} R_{\infty}^n$ and $C_{\infty}^n
\otimes_h E = C_{\infty}^n \otimes_{\mbox{\scriptsize min}} E$ and
the complete contractions $$\begin{array}{l} R_{\infty}^n
\otimes_h (E \otimes_{\mbox{\scriptsize min}} R_{\infty}^n)
\longrightarrow R_{\infty}^n \otimes_{\mbox{\scriptsize min}} (E
\otimes_{\mbox{\scriptsize min}} R_{\infty}^n), \\ (C_{\infty}^n
\otimes_{\mbox{\scriptsize min}} E) \otimes_h C_{\infty}^n
\longrightarrow (C_{\infty}^n \otimes_{\mbox{\scriptsize min}} E)
\otimes_{\mbox{\scriptsize min}} C_{\infty}^n. \end{array}$$
Hence, it suffices to show our claim for the diagonals of
$R_{\infty}^{n^2} \otimes_{\mbox{{\scriptsize min}}} E$ and
$C_{\infty}^{n^2} \otimes_{\mbox{{\scriptsize min}}} E$. In the
first case the diagonal is $R_{\infty}^n
\otimes_{\mbox{{\scriptsize min}}} E$ while in the second case is
$C_{\infty}^n \otimes_{\mbox{{\scriptsize min}}} E$. By the
injectivity of the minimal tensor product and since
$\ell_{\infty}^n$ carries the minimal operator space structure,
the assertion follows. This completes the proof. \fin

The following result can be regarded as the dual version of
Proposition \ref{Proposition-Lambda-Cap}, where the spaces
$J_{p,q}^n(M_l;E)$ are replaced by the spaces $K_{p,q}^n (M_l;E)$.
Here we skip the assumption that $q=1$ and we work in the range $1
\le p \le q \le \infty$.

\begin{proposition} \label{Proposition-Lambda-Sum}
Let $E$ be an operator space and let $1 \le p \le q \le \infty$.
Then, the following map is a complete contraction $$\Lambda_{pq}:
\sum_{k=1}^n \delta_k \otimes x_k \in K_{p,q}^n (M_l;E)
\longmapsto \sum_{k=1}^n \delta_k \otimes \pi_k(x_k) \in
L_p(\tau_{l^n}; \ell_q^n(E)).$$
\end{proposition}

\dem Let $t = l^{\frac{1}{2p} - \frac{1}{2q}}$, regarding again
$K_t(C_p^{nl},C_q^{nl}) \otimes_h E \otimes_h
K_t(R_p^{nl},R_q^{nl})$ as a space of $n \times n$ matrices with
entries in $M_l \otimes E$, we consider its diagonal subspace
$\mathcal{K}_{p,q}^n(M_l;E)$. By Remark \ref{Remark-K-Diagonal},
it suffices to check that the mapping $$\widetilde{\Lambda}_{pq}:
\sum_{k=1}^n e_{kk} \otimes x_k \in \mathcal{K}_{p,q}^n (M_l;E)
\longmapsto \sum_{k=1}^n \delta_k \otimes \pi_k(x_k) \in
L_p(\tau_{l^n}; \ell_q^n(E)),$$ satisfies
$\|\widetilde{\Lambda}_{pq}\|_{cb} \le \lambda_p$. Since the
diagonal projection $\mathrm{P}$ is a complete contraction, it
suffices to prove this estimate for the diagonal in each of the
following spaces $$C_p^{nl} \ten_h E \ten_h R_p^{nl} \pl , \pl t
C_p^{nl} \ten_h E \ten_h R_q^{nl} \pl , \pl t C_q^{nl} \ten_h E
\ten_h R_p^{nl} \pl ,\pl t^2 C_q^{nl} \ten_h E \ten_h R_q^{nl}.$$
Notice that, given a scalar $\gamma$ and an operator space $F$, we
denote by $\gamma F$ the operator space with operator space
structure given by $$\|f\|_{M_m \otimes_{\mbox{\scriptsize min}}
\gamma F} = \gamma \|f\|_{M_m \otimes_{\mbox{\scriptsize min}}
F}.$$ The first one is $$\big\| \widetilde{\Lambda}_{pq}:
\ell_p^n(S_p^l(E)) \rightarrow L_p(\tau_{l^n}; \ell_q^n(E))
\big\|_{cb} \le \lambda_p.$$ This estimate obviously holds for
$p=q$ and, since the identity map $\ell_p^n(E) \rightarrow
\ell_q^n(E)$ is a complete contraction, the desired estimate
follows. For the last one, we note that $$\big\|
\widetilde{\Lambda}_{pq}: t^2 \ell_q^n(S_q^l(E)) \rightarrow
L_q(\tau_{l^n}; \ell_q^n(E)) \big\|_{cb} \le \lambda_p.$$
Moreover, since we are using a probability measure, we know that
the identity map $L_q(\tau_m; F) \rightarrow L_p(\tau_m; F)$ is a
complete contraction. Therefore, the desired estimate for the last
case holds. For the second and third terms, we use a similar
trick. We claim that the identity mappings
\begin{eqnarray*}
L_{(2p,2q)}(\tau_m; F_1) & \longrightarrow & L_p(\tau_m; F_1)
\\ L_{(2q,2p)}(\tau_m; F_2) & \longrightarrow &
L_p(\tau_m; F_2)
\end{eqnarray*}
are complete contractions. Namely, by complex interpolation it
reduces to the case $p=1$ and $q = \infty$. However, if we rescale
these mappings to replace $\tau_m$ by $\mbox{tr}_m$, this case
follows easily by the injectivity of the Haagerup tensor product
and the well-know estimates $$\big\| id: R_{\infty}^m \to
C_{\infty}^m \big\|_{cb} \le \sqrt{m} \quad \mbox{and} \quad
\big\| id: C_{\infty}^m \to R_{\infty}^m \big\|_{cb} \le
\sqrt{m}.$$ We take $m=l^n$ and the operator spaces $F_1 =
\ell_{(2p,2q)}^n (E)$ and $F_2 = \ell_{(2q,2p)}^n (E)$. According
to Lemma \ref{Lemma-Contractions}, it suffices to prove the
following estimates
\begin{eqnarray*}
\big\| \widetilde{\Lambda}_{pq}: t \, \ell_{(2p,2q)}^n \big(
S_{(2p,2q)}^l(E) \big) \rightarrow L_{(2p,2q)}(\tau_{l^n};
\ell_{(2p,2q)}^n(E)) \big\|_{cb} & \le & \lambda_p \\ \big\|
\widetilde{\Lambda}_{pq}: t \, \ell_{(2q,2p)}^n \big(
S_{(2q,2p)}^l(E) \big) \rightarrow L_{(2q,2p)}(\tau_{l^n};
\ell_{(2q,2p)}^n(E)) \big\|_{cb} & \le & \lambda_p.
\end{eqnarray*}
Since both estimates can be treated in a similar way, we just
prove the first one. Given a positive integer $m$, let us consider
a diagonal matrix $$x = \sum_{k=1}^n e_{kk} \otimes x_k \in M_n
\otimes S_{(2p,2q)}^{ml}(E).$$ According to Lemma
\ref{Lemma-Asymmetric} (iii), the following identities hold for
$\frac{1}{r} = \frac{1}{2p} + \frac{1}{2q}$
\begin{eqnarray*}
\big\| \widetilde{\Lambda}_{pq}(x) \big\|_{(2p,2q)} & = & \Big(
\sum_{k=1}^n \|\pi_k(x_k)\|_{L_{(2p,2q)}(\mathrm{tr}_m \otimes
\tau_{l^n};E)}^r \Big)^{1/r} \\ & = & \Big( \sum_{k=1}^n
\|x_k\|_{L_{(2p,2q)}(\mathrm{tr}_m \otimes \tau_l; E)}^r
\Big)^{1/r}
\\ & = & \lambda_p \, t \|x\|_{\ell_{(2p,2q)}^n
(S_{(2p,2q)}^{ml}(E))},
\end{eqnarray*}
where $\| \ \|_{(2p,2q)}$ denotes the norm on the space
$L_{(2p,2q)}(\mbox{tr}_m \otimes \tau_{l^n};
\ell_{(2p,2q)}^n(E))$. Thus, applying the second identity of Lemma
\ref{Lemma-Asymmetric} (ii), the assertion follows. \fin

Once we have seen the estimates for intersections and sums given
in Propositions \ref{Proposition-Lambda-Cap} and
\ref{Proposition-Lambda-Sum}, we are in position to prove the
complete equivalence of norms (\ref{Equation-S(p,1)}) for sums of
independent noncommutative random variables in $L_p(\ell_1(E))$.

\begin{theorem} \label{Theorem-Main(q=1)}
Let $E$ be an operator space and let $1 \le p < \infty$. Then, the
map $$\Lambda_{p1}:  \sum_{k=1}^n \delta_k \otimes x_k \in
J_{p,1}^n(M_l;E) \longmapsto \sum_{k=1}^n \delta_k \otimes
\pi_k(x_k) \in L_p(\tau_{l^n}; \ell_1^n(E))$$ is a complete
isomorphism onto a completely complemented subspace. Similarly,
the same holds for the map $$\Lambda_{p'\infty}: \sum_{k=1}^n
\delta_k \otimes x_k \in K_{p',\infty}^n(M_l; E) \longmapsto
\sum_{k=1}^n \delta_k \otimes \pi_k(x_k) \in L_{p'}(\tau_{l^n};
\ell_\infty^n(E))$$ for $1 < p' \le \infty$. Moreover, the
$cb$-distance constants are independent of $l$ and $n$.
\end{theorem}

\dem It is clear that we can assume $E$ to be a finite-dimensional
operator space. In particular, all the spaces we shall consider
along the proof will be of finite dimension and hence reflexive.
Now the duality theory for the Haagerup tensor product, see for
instance the Chapter 5 of \cite{P3}, provides a complete isometry
$$\mathrm{S}: \big( J_t(C_p^{nl}, C_1^{nl}) \otimes_h E \otimes_h
J_t(R_p^{nl}, R_1^{nl}) \big)^* \rightarrow
K_{t^{-1}}(C_{p'}^{nl}, C_{\infty}^{nl}) \otimes_h E^* \otimes_h
K_{t^{-1}}(R_{p'}^{nl}, R_{\infty}^{nl}).$$ On the other hand,
according to Lemma \ref{Lemma-Projection} and Remark
\ref{Remark-K-Diagonal}, the projection onto the diagonal is
always a complete contraction. Therefore, we obtain the following
completely isometric isomorphism $$J_{p,1}^n(M_l; E)^* =
K_{p',\infty}^n(M_l; E^*).$$ Indeed, if $\mathrm{P}$ denotes the
diagonal projection and $\mathrm{T} = u^{-1} \circ \mathrm{P}$
where $u$ stands for the linear mapping considered in Lemma
\ref{Lemma-Projection}, then the mapping $$\mathrm{T} \circ
\mathrm{S} \circ \mathrm{T}^*: J_{p,1}^n(M_l; E)^* \longrightarrow
K_{p',\infty}^n(M_l; E^*)$$ is a completely isometric isomorphism.
Here, the duality is given by $$\langle a,b \rangle = \left\langle
\sum_{k=1}^n \delta_k \otimes (a_k \otimes e_k), \sum_{k=1}^n
\delta_k \otimes (b_k \otimes e_k^*) \right\rangle = \sum_{k=1}^n
\tau_l(a_k^{\mbox{t}} b_k^{}) \langle e_k, e_k^* \rangle.$$ Thus,
we obviously have $$\langle \Lambda_{p1} (a), \Lambda_{p'\infty}
(b) \rangle = \langle a,b \rangle \qquad \forall \ a \in
J_{p,1}^n(M_l;E), \ b \in K_{p',\infty}^n(M_l;E^*).$$
Consequently, the map $\Lambda_{p'\infty}^* \Lambda_{p1}^{}$ is
the identity on $J_{p,1}^n(M_l;E)$. In particular, by Propositions
\ref{Proposition-Lambda-Cap} and \ref{Proposition-Lambda-Sum},
$\Lambda_{p1}$ becomes a complete isomorphism with constants not
depending on the dimensions. Moreover, its image is a completely
complemented subspace since $\Lambda_{p1}^{} \Lambda_{p'\infty}^*$
is a completely bounded projection with $\|\Lambda_{p1}^{}
\Lambda_{p'\infty}^*\|_{cb} \le cp$. This proves the assertions
for $\Lambda_{p1}$, but the arguments for $\Lambda_{p'\infty}$ are
similar. \fin

\begin{remark} \label{Remark-Emphasis}
\emph{Let us state Theorem \ref{Theorem-Main(q=1)} in a more
explicit way. To that aim, we introduce some notation. If
$\frac{1}{\gamma_{rs}} = \frac1r + \frac1s$, we define
$$\begin{array}{l} \displaystyle \|x\|_{p,q}^{\cap} = \max_{r,s
\in \{2p,2q\}} \left\{ \Big( \sum_{k=1}^n
\|x_k\|_{L_{(r,s)}(\tau_l;E)}^{\gamma_{rs}} \Big)^{1/\gamma_{rs}}
\right\}, \\ \displaystyle \|x\|_{p,q}^{\Sigma} =
\inf_{{\mbox{\tiny $\displaystyle x = \sum_{r,s} x_{rs}$}}}
\left\{ \sum_{r,s} \Big( \sum_{k=1}^n
\|x_{rs}^k\|_{L_{(r,s)}(\tau_l;E)}^{\gamma_{rs}}
\Big)^{1/\gamma_{rs}} \, \Big| \ r,s \in \{2p,2q\} \right\}.
\end{array}$$
Then, recalling the meaning of $\le_{cb}$ from the Introduction,
we have
\begin{itemize}
\item Given $1 \le p < \infty$, we have $$\|x\|_{p,1}^{\cap}
\le_{cb} \Big\| \sum_{k=1}^n \delta_k \otimes \pi_k(x_k)
\Big\|_{L_p(\ell_1^n(E))} \le_{cb} cp \ \|x\|_{p,1}^{\cap}.$$
\item Given $1 < p' \le \infty$, we have
$$\frac{1}{cp} \ \|x\|_{p',\infty}^{\Sigma} \le_{cb} \Big\|
\sum_{k=1}^n \delta_k \otimes \pi_k(x_k)
\Big\|_{L_{p'}(\ell_{\infty}^n(E))} \le_{cb}
\|x\|_{p',\infty}^{\Sigma}.$$
\end{itemize}}
\end{remark}

\section{A $cb$ embedding of $S_q^n$ into $S_p(\ell_q^m)$}
\label{Section4}

We begin by stating a complementation result for the subspace of
$J_{p,q}^n(M_l;E)$ given by constant diagonal matrices. As we
shall see immediately, this result plays a relevant role in the
embeddings we want to consider.

\begin{lemma} \label{Lemma-Constant-Diagonal}
Let $1\le p,q\le \infty$ and let $t=\kla \frac{n}{l}
\mer^{\frac{1}{2p}-\frac{1}{2q}}$ with $l$ and $n$ positive
integers. Then, the map $$\mathrm{T}: J_t(C_q^l,C_p^l) \ten_h E
\ten_h J_t(R_q^l,R_p^l) \longrightarrow J_{p,q}^n(M_l;E)$$ defined
by $$\mathrm{T}(x) = \kla \frac{n}{l}\mer^{-1/q} \Big(
\sum_{k=1}^n e_{kk} \otimes x \Big)$$ is a complete isometry. The
image of $\mathrm{T}$ is completely contractively complemented.
\end{lemma}

\dem To see that the image of $\mathrm{T}$ is completely
contractively complemented in $J_{p,q}^n(M_l;E)$, we consider the
following projection $$\mathrm{P} (x_1, x_2, \ldots, x_n) = \Big(
\frac{1}{n} \sum_{k=1}^n x_k, \frac{1}{n} \sum_{k=1}^n x_k,
\ldots, \frac{1}{n} \summ_{k=1}^n x_k \Big).$$ Then, it suffices
to see that $\mathrm{P}$ is a complete contraction in
$$\ell_{(r,s)}^n (L_{(r,s)} (\tau_l; E))$$ whenever $r,s \in
\{2p,2q\}$. It is clear that, given any operator space $E$, the
projection $\mathrm{P}$ is contractive in these four spaces. Then,
the complete contractivity follows easily from Lemma
\ref{Lemma-Asymmetric} (ii) and the obvious Fubini type results.
Now, given $r,s \in \{2p,2q\}$, let $\xi_{rs} = \delta_{r,2p} +
\delta_{s,2p}$. To see that $\mathrm{T}$ is a complete isometry,
it suffices to check that $$\mathrm{T}: t^{\xi_{rs}} S_{(r,s)}^l
(E) \longrightarrow \ell_{(r,s)}^n (L_{(r,s)} (\tau_l; E))$$ is a
complete isometry for any $r,s \in \{2p,2q\}$. However, this
follows one more time as a consequence of Lemma
\ref{Lemma-Asymmetric} (ii) and (iii). \fin

The following theorem provides an embedding of the Schatten class
$S_q^n(E)$ into $L_p(\EME, \tau; \ell_q^m(E))$ with uniformly
bounded $cb$-distance constants.

\begin{theorem} \label{Theorem-Embedding}
Let $1 \le q \le p < \infty$. Then, given any positive integer $n
\ge 1$ and any operator space $E$, the following mapping is a
complete isomorphism onto a completely complemented subspace
$$\Phi_{pq}: x \in S_q^n(E) \longmapsto \frac{1}{n^{1/q}}
\sum_{k=1}^{n^2} \delta_k \otimes \pi_k(x) \in L_p(\tau_{n^{n^2}};
\ell_q^{n^2}(E)).$$ Moreover, $\|\Phi_{pq}\|_{cb} \le c p$ while
the inverse mapping $\Phi_{pq}^{-1}$ is completely contractive.
\end{theorem}

\dem By Lemma \ref{Lemma-CB(Cr,Cs)}, $J_t(C_1^n,C_p^n) =
R_{\infty}^n$ and $J_t(R_1^n,R_p^n) = C_{\infty}^n$ for $t =
n^{\frac{1}{2p}-\frac12}$. In particular, we can write $$S_1^n(E)
= J_t(C_1^n,C_p^n) \ten_h E \ten_h J_t(R_1^n,R_p^n).$$ Then,
Proposition \ref{Proposition-Lambda-Cap} and Lemma
\ref{Lemma-Constant-Diagonal} give that $$\Phi_{p1}: S_1^n(E)
\longrightarrow L_p(\tau_{n^{n^2}}; \ell_1^{n^2}(E))$$ is a $cb$
embedding with $\|\Phi_{p1}\|_{cb} \le cp$. That is, the upper
estimate holds for $q=1$. On the other hand, the map $$\Phi_{pp}:
S_p^n(E) \longrightarrow L_p(\tau_{n^{n^2}}; \ell_p^{n^2}(E))$$ is
clearly a complete isometry. Hence, for the general case, the
upper estimate follows by complex interpolation. In order to see
that the image of the mapping $\Phi_{pq}$ is completely
complemented and $\Phi_{pq}^{-1}$ is completely contractive, we
observe again that, by elementary properties of the local theory,
we have $$S_{q'}^n(E^*) = \big( J_t(C_q^n,C_p^n) \ten_h E \ten_h
J_t(R_q^n,R_p^n) \big)^*$$ for $t =
n^{\frac{1}{2p}-\frac{1}{2q}}$. Thus, if
$\mathcal{C}_{p',q'}^{n^2} (M_n; E^*)$ stands for the subspace of
$K_{p',q'}^{n^2}(M_n; E^*)$ of constant diagonals, Lemma
\ref{Lemma-Constant-Diagonal} and duality give $$S_{q'}^n(E^*) =
\frac{1}{n^{1/q}} \, \mathcal{C}_{p',q'}^{n^2} (M_n; E^*).$$ In
particular, Proposition \ref{Proposition-Lambda-Sum} gives that
$$\Phi_{p'q'}: S_{q'}^n(E^*) \longrightarrow
L_{p'}(\tau_{n^{n^2}}; \ell_{q'}^{n^2}(E^*))$$ is completely
contractive. Finally, we observe that $$\big\langle
\Phi_{pq}(a\ten e),\Phi_{p'q'}(b\ten e^*) \big\rangle = \frac1n
\sum_{k=1}^{n^2} \tau_n(a^{t}b) \langle e,e^*\rangle = \langle a
\ten e,b \ten e^* \rangle.$$ Hence, since $\Phi_{p'q'}^*
\Phi_{pq}$ is the identity and $\Phi_{pq} \Phi_{p'q'}^*$ is a
projection, we are done. \fin

\begin{remark} \label{Remark-Embedding-Duality}
\emph{By simple dual arguments, it is not difficult to check that
Theorem \ref{Theorem-Embedding} holds for $1 < p \le q \le
\infty$, with $\Phi_{pq}$ completely contractive and
$\|\Phi_{pq}^{-1}\|_{cb} \le cp$. Namely, we first recall that
$$S_{\infty}^n(E) = K_t(C_{\infty}^n,C_p^n) \otimes_h E \otimes_h
K_t(R_{\infty}^n,R_p^n) \qquad \mbox{for} \quad t=
n^{\frac{1}{2p}}.$$ Then, by Theorem \ref{Theorem-Main(q=1)} and
the dual version of Lemma \ref{Lemma-Constant-Diagonal} for the
$K$ functional, the complete contractivity of $\Phi_{p\infty}$
holds. Finally, we end by interpolation and duality.}
\end{remark}

\begin{remark} \label{Remark-Optimal}
\emph{Rescaling Theorem \ref{Theorem-Embedding}, we get an
embedding $\Psi_{pq}: S_q^n \rightarrow S_p(\ell_q^{m})$. In fact,
we have taken $m$ to be $n^2$. As we shall see in Section
\ref{Section6}, when seeking for $cb$-embeddings with uniformly
bounded constants, the choice $m = n^2$ is optimal.}
\end{remark}

\section{$\mathrm{K}$-convex operator spaces}
\label{Section5}

The theory of type and cotype is essential to study some geometric
properties of Banach spaces. The operator space analog of that
theory has been recently initiated in some works summarized in
\cite{Pa2}. The aim of this section is to explore the relation
between $\mathrm{B}$-convexity and $\mathrm{K}$-convexity in the
category of operator spaces.

\subsection{A variant of the embedding theorem}

In this paragraph, we study the inverse of $\Phi_{pq}$ when we
impose on $\ell_q$ its minimal operator space structure. The
resulting mapping will the key in the operator space analog of
Pisier's equivalence between $\mathrm{B}$-convex and
$\mathrm{K}$-convex spaces.

\begin{lemma} \label{Lemma-Antilinear}
The following map extends to an anti-linear isometry $$\mathrm{T}:
\sum_{k=1}^n a_k \otimes e_k \in L_p(\tau_n;\min(E)) \longmapsto
\sum_{k=1}^n a_k^* \otimes \overline{e}_k \in
L_p(\tau_n;\overline{\min(E)}).$$ Here, $\overline{\min(E)}$
stands for the complex conjugate operator space as defined in
\textnormal{\cite{P3}}.
\end{lemma}

\dem Since $\min (E)$ embeds completely isometrically in
$\ell_\infty$, we take $E$ to be $\ell_\infty$. Under this
assumption, the result is clear for $p=\infty$. Namely, given $x =
(x_n)_{n \ge 1}$ in $L_\infty(\tau_n; \ell_\infty)$, we have
$$\|\mathrm{T}(x)\| = \sup_{n \ge 1} \|x_n^*\| = \|x\|.$$ Now, if
$x \in L_p(\tau_n; \ell_\infty)$, there exist $a,b \in
L_{2p}(\tau_n)$ and $y \in L_{\infty}(\tau_n; \ell_{\infty})$ such
that $x = ayb$ and $$\|a\|_{2p} \|y\|_{\infty} \|b\|_{2p} < (1 +
\eps) \|x\|.$$ Therefore $$\|\mathrm{T}(x)\| \le \|b^*\|_{2p}
\|y\|_{\infty} \|a^*\|_{2p} < (1 + \eps) \|x\|.$$ Since $\eps
> 0$ is arbitrary, the assertion follows easily. This completes
the proof. \fin

Let us consider the operator space $\widetilde{\mathsf{F}}_{pq}^n$
defined as the image of $S_q^n$ under the map $\Phi_{pq}$, with
the operator space structure inherited from $$L_p(\tau_{n^{n^2}};
\min (\ell_q^{n^2})).$$

\begin{proposition} \label{Proposition-Min}
The estimate
$\|\Phi_{pq}^{-1}\|_{\mathcal{B}(\widetilde{\mathsf{F}}_{pq}^n,S_q^n)}
\le 2$ holds for any $n \ge 1$.
\end{proposition}

\dem We first consider a self-adjoint matrix $x$ in $S_q^n$. Then
taking $m=n^2$, the sequence $\pi_1(x), \pi_2(x), \ldots,
\pi_m(x)$ lies in a commutative subalgebra of $M_{n^m}$. In fact,
using the spectral theorem, we can write $x = u^*
\mathrm{d}_{\lambda} u$ where $\mathrm{d}_{\lambda}$ stands for
the matrix of eigenvalues of $x$ and $u$ is unitary. In
particular, after multiplication by $u^{\ten m}$ from the left and
by $(u^*)^{\ten m}$ from the right, we may assume that
$$\sum_{k=1}^m \delta_k \otimes \pi_k(x)$$ is a diagonal matrix.
In that case, we may apply Corollary 1.3 of \cite{P2} to obtain
$$\Big\|\sum_{k=1}^m \delta_k \ten \pi_k(x)
\Big\|_{L_p(\tau_{n^m}; \min(\ell_q^m))} = \Big\| \sum_{k=1}^m
\delta_k \ten \pi_k(x) \Big\|_{L_p(\tau_{n^m}; \ell_q^m)}.$$
Therefore, Theorem \ref{Theorem-Embedding} gives $$\|x\|_{S_q^n}
\le
\|\Phi_{pq}(x)\|_{\widetilde{\mathsf{F}}_{pq}^n}.$$ For arbitrary
$x$, we consider its decomposition into self-adjoint elements $$a
= \frac12 (x+x^*) \quad \mbox{and} \quad b = \frac{1}{2i}
(x-x^*).$$ Then, we deduce from Lemma \ref{Lemma-Antilinear} that
$$\|\Phi_{pq}(a)\|_{\widetilde{\mathsf{F}}_{pq}^n} \le \frac12
\|\Phi_{pq}(x)\|_{\widetilde{\mathsf{F}}_{pq}^n} + \frac12
\|\Phi_{pq}(x)^*\|_{\widetilde{\mathsf{F}}_{pq}^n} \le
\|\Phi_{pq}(x)\|_{\widetilde{\mathsf{F}}_{pq}^n}.$$ Obviously, the
same estimate holds for $b$. Thus, we obtain the desired estimate.
\fin

Let us consider an infinite dimensional operator space $E$ and a
family of finite dimensional operator spaces $\mathcal{A} = \big\{
A_n \, \big| \ n \ge 1 \big\}$. We shall say that the family
$\mathcal{A}$ \emph{embeds semi-completely uniformly} in $E$, and
we shall write $\mathcal{A} \prec E$, when there exists a constant
$\mathrm{c}$ and embeddings $\Lambda_n: A_n \rightarrow E$ such
that $$\|\Lambda_n\|_{cb} \|\Lambda_n^{-1}\| \le \mathrm{c} \qquad
\mbox{for all} \qquad n \ge 1.$$

\begin{corollary} \label{Corollary-Semi-Embedding}
Let $1 < p < \infty$ and $1 \le q \le \infty$. Then, we have
$$\Big\{ \ell_q^n \, \Big| \ n \ge 1 \Big\} \prec S_p(E)
\Rightarrow \Big\{ S_q^n \, \Big| \ n \ge 1 \Big\} \prec S_p(E).$$
\end{corollary}

\dem By hypothesis, there exist $\mathrm{c}_1 > 1$ and embeddings
$\Lambda_n: \ell_q^n \rightarrow S_p(E)$ such that
$\|\Lambda_n\|_{cb} \|\Lambda_n^{-1}\| \le \mathrm{c}_1$ for each
positive integer $n$. Let $F_n$ denote the image
$\Lambda_n(\ell_q^n)$ of $\Lambda_n$ in $S_p(E)$. On the other
hand, according to Theorem \ref{Theorem-Embedding}, we know how to
construct linear isomorphisms $$\Phi_n: S_q^n \rightarrow
\mathsf{F}_{pq}^n \subset S_p(\ell_q^{n^2}) \qquad \mbox{such
that} \qquad \|\Phi_n\|_{cb} \le \mathrm{c}_2$$ for some constant
$\mathrm{c}_2$ independent of $n$. Moreover, let
$\widetilde{\mathsf{F}}_{pq}^n$ be the image of $\Phi_n$ endowed
with the operator space structure inherited from $$S_p(\min
(\ell_q^{n^2})).$$ Then, if $\Psi_n: \widetilde{\mathsf{F}}_{pq}^n
\to S_q^n$ stands for $\Phi_n^{-1}$, Proposition
\ref{Proposition-Min} gives $\|\Psi_n\| \le \mathrm{c}_3$ for some
constant $\mathrm{c}_3$ independent of $n$. Let us define
$$\widetilde{\Lambda}_n: S_q^n \rightarrow S_p(E) \qquad \mbox{by}
\qquad \widetilde{\Lambda}_n = \big( id \otimes \Lambda_{n^2}
\big) \circ \Phi_n.$$ Then we have
\begin{eqnarray*}
\|\widetilde{\Lambda}_n\|_{cb} \|\widetilde{\Lambda}_n^{-1}\| &
\le & \|\Lambda_{n^2}\|_{cb} \|\Phi_n\|_{cb} \|\Psi_n\|
\|\Lambda_{n^2}^{-1}\|_{\mathcal{CB}(F_{n^2}, \min
(\ell_q^{n^2}))}
\\ & = & \|\Lambda_{n^2}\|_{cb} \|\Phi_n\|_{cb} \|\Psi_n\|
\|\Lambda_{n^2}^{-1}\| \le \mathrm{c}_1 \mathrm{c}_2 \mathrm{c}_3.
\end{eqnarray*}
Since the constant $\mathrm{c}_1 \mathrm{c}_2 \mathrm{c}_3$ does
not depend on $n$, the assertion follows. \fin

\subsection{$\mathrm{OB}$-convexity and $\mathrm{OK}$-convexity}

Let us start by defining the notion of $\mathrm{OB}$-convex
operator space. Following \cite{Pa}, let us fix a family
$\mathbf{d}_{\Sig} = \big\{ \Degree: \sig \in \Sig \big\}$ of
positive integers indexed by an infinite set $\Sig$ and, given a
finite subset $\Gamma$ of $\Sig$, let $$\Delta_{\Gamma} =
\sum_{\sig \in \Gamma} \Degree^2.$$ An operator space $E$ is
called \textbf{$\mathrm{OB}_{\Sig}$-convex} if there exists a
finite subset $\Gamma$ of $\Sig$ and certain $0 < \delta \le 1$
such that, for any family $$\Big\{ A^{\sig} \in M_{\Degree}
\otimes S_2(E) \Big\}_{\sig \in \Gamma}^{\null},$$ we have
$$\frac{1}{\Delta_{\Gamma}} \inf_{B^{\sig} unitary} \Big\|
\sum_{\sig \in \Gamma} \Degree \mbox{tr}(A^{\sig} B^{\sig})
\Big\|_{S_2(E)}^{\null} \le (1 - \delta) \, \max_{\sig \in
\Gamma}\|A^{\sig}\|_{M_{\Degree}(S_2(E))}.$$

If we replace the Schatten class $S_2(E)$ above by $S_p(E)$ we get
an equivalent notion whenever $1 < p < \infty$, see \cite{Pa}.
This definition is inspired by Beck's original notion for Banach
spaces, which corresponds to the \emph{commutative set of
parameters} $\Sig_0 = \N$ with $d_{\sigma} = 1$ for all $\sig \in
\Sig_0$. Our definition depends a priori on the set of parameters
$(\Sig, \mathbf{d}_{\Sig})$. However, we shall see below that
there is no dependence on $\Sig$. On the other hand, we also need
to provide an operator space analog of the property of containing
(uniformly) finite dimensional $L_1$ spaces. However, this time we
need to allow the noncommutative $L_1$'s to appear in the
definition. Given an operator space $E$, a set of parameters
$(\Sig, \mathbf{d}_{\Sig})$ and $1 \le p < \infty$, we define the
spaces $$\Lebesgue_p(\Sig;E) = \left\{ A \in \prod_{\sig \in \Sig}
M_{\Degree} \ten E: \Big( \sum_{\sig \in \Sig} \Degree
\|A^{\sig}\|_{S_p^{\Degree}(E)}^p \Big)^{1/p} < \infty \right\}.$$
We impose on $\Lebesgue_p(\Sig;E)$ its natural operator space
structure, see Chapter 2 of \cite{P2} for the details. We shall
write $\Lebesgue_p(\Sig)$ for the scalar-valued case. We shall say
that $S_p(E)$ \textbf{contains $\Lebesgue_1(\Gamma)$'s
semi-completely $\lambda$-uniformly} if, for each finite subset
$\Gamma$ of $\Sig$, there exists a linear embedding
$\Lambda_{\Gamma}: \Lebesgue_1(\Gamma) \rightarrow S_p(E)$ such
that $$\|\Lambda_{\Gamma}\|_{cb} \|\Lambda_{\Gamma}^{-1}\| \le
\lambda.$$ In other words, if $$\Big\{ \Lebesgue_1(\Gamma) \,
\Big| \ \Gamma \ \mbox{finite} \Big\} \prec S_p(E).$$ The
following is the analog of a well-known result for Banach spaces,
see \cite{Pa}.

\begin{remark} \label{Remark-Uniform-lambda}
\emph{Given an operator space $E$, the following are equivalent:
\begin{itemize}
\item[i)] $S_p(E)$ contains $\Lebesgue_1(\Gamma)$'s semi-completely
$\lambda$-uniformly for any $\lambda > 1$.
\item[ii)] $S_p(E)$ contains $\Lebesgue_1(\Gamma)$'s semi-completely
$\lambda$-uniformly for some $\lambda > 1$.
\end{itemize}}
\end{remark}

Finally we recall, as have already done in the Introduction, that
an operator space $E$ will be considered
\textbf{$\mathrm{OK}$-convex} whenever the vector-valued Schatten
class $S_2(E)$ is $\mathrm{K}$-convex when regarded as a Banach
space.

\begin{remark} \label{Remark-p-Independence}
\emph{The given definition of $\mathrm{OK}$-convexity is a bit
more flexible. Indeed, an operator space $E$ is
$\mathrm{OK}$-convex if and only if $S_p(E)$ is a
$\mathrm{K}$-convex Banach space for some (any) $1 < p < \infty$.
This follows from the fact that, given $1 < p < \infty$, the
Schatten class $S_p(E)$ is $\mathrm{K}$-convex if and only
$S_2(E)$ is $\mathrm{K}$-convex. Indeed, it follows from
\cite{Pi,P0} that Banach space $\mathrm{K}$-convexity is stable by
complex interpolation assuming only that one of the endpoint
spaces is $\mathrm{K}$-convex. Now assume that $S_2(E)$ is
$\mathrm{K}$-convex and let $1 < p < \infty$. If $p < 2$ (resp. $p
> 2$) we have $$S_p(E) = [S_2(E), S_1(E)]_{\theta} \qquad
\big(\mbox{resp.} \ S_p(E) = [S_2(E), S_{\infty}(E)]_{\theta}
\big)$$ for some $0 < \theta < 1$. Therefore, we find by complex
interpolation that $S_p(E)$ is also a $\mathrm{K}$-convex Banach
space. A similar argument shows that $S_2(E)$ is a
$\mathrm{K}$-convex Banach space whenever $S_p(E)$ is also
$\mathrm{K}$-convex. Thus our claim follows.}
\end{remark}

\begin{remark} \label{Remark-K-Sig-Covexity}
\emph{In \cite{Pa} it was given an a priori more general notion of
$\mathrm{K}$-convexity for operator spaces. Namely, let $(\Omega,
\mathsf{A}, \mu)$ be a probability space with no atoms. Then,
following \cite{MP} we define the \emph{quantized Gauss system}
associated to $(\Sig, \mathbf{d}_{\Sig})$ as a collection of
matrix-valued functions $$\mathbf{G}_{\Sig} = \Big\{
\gamma^{\sig}: \Omega \rightarrow M_{\Degree} \Big\}_{\sig \in
\Sig} \quad \mbox{where} \quad \gamma^{\sig} =
\frac{1}{\sqrt{\Degree}} \,\ \Big(\,\ \gau_{ij}^{\sig} \,\
\Big).$$ Here, the functions $\gau_{ij}^{\sig}: \Omega \rightarrow
\mathbb{C}$ form a family, indexed by $1 \le i,j \le \Degree$ and
$\sig \in \Sig$, of independent standard complex-valued gaussian
random variables. Given a function $f \in L_2(\Omega;E)$, we can
consider the Fourier coefficients of $f$ with respect to the
quantized Gauss system $$\widehat{f}_{\mathbf{G}}(\sig) =
\int_{\Omega} f(\omega) \gamma^{\sig}(\omega)^{\ast} d
\mu(\omega).$$ This gives rise to the Gauss projection defined
below $$\mathrm{P}_\mathbf{G}: f \in L_2(\Omega;E) \longmapsto
\sum_{\sig \in \Sig} \Degree
\mbox{tr}(\widehat{f}_{\mathbf{G}}(\sig) \gamma^{\sig}) \in
L_2(\Omega;E).$$ An operator space $E$ is called
\emph{$\mathrm{OK}_{\Sig}$-convex} if the Gauss projection
associated to the parameters $(\Sig, \mathbf{d}_{\Sig})$ is a
completely bounded map. However, recalling the definition of the
quantized Gauss system, we can write $$\sum_{\sig \in \Sig}
\Degree \mbox{tr}(\widehat{f}_{\mathbf{G}}(\sig) \gamma^{\sig}) =
\sum_{\sig \in \Sig} \sum_{i,j = 1}^{\Degree} \int_{\Omega}
f(\omega) \overline{\gau_{ij}^{\sig}(\omega)} d \mu(\omega) \,
\gau_{ij}^{\sig}.$$ Therefore, since now the right hand side can
be regarded as the classical Gauss projection, it turns out that
the notion of $\mathrm{OK}_{\Sig}$-convexity does not depend on
the set $(\Sig, \mathbf{d}_{\Sig})$, so that we shall simply use
in the sequel the term $\mathrm{OK}$-convex, without any explicit
reference to the set of parameters $(\Sig, \mathbf{d}_{\Sig})$.}
\end{remark}

\begin{remark}
\emph{We can replace $L_2(\Omega;E)$ above by $L_p(\Omega;E)$ for
any $1 < p < \infty$.}
\end{remark}

\begin{theorem} \label{Theorem-Pisier}
Given an operator space $E$, the following are equivalent:
\begin{itemize}
\item[i)] $E$ is $\mathrm{OK}$-convex.
\item[ii)] $E$ is $\mathrm{OB}_{\Sig}$-convex for some $($any$)$
set of parameters $(\Sig, \mathbf{d}_{\Sig})$.
\item[iii)] $S_p(E)$ does not
contain $\ell_1^n$'s uniformly for some $($any$)$ $1 < p <
\infty$.
\item[iv)] $S_p(E)$ does not contain $\Lebesgue_1(\Gamma)$'s
semi-completely for some $($any$)$ $1 < p < \infty$.
\end{itemize}
\end{theorem}

\dem By definition, $E$ is $\mathrm{OK}$-convex if and only if
$S_p(E)$ is a $\mathrm{K}$-convex Banach space for some (any) $1 <
p < \infty$, see Remark \ref{Remark-p-Independence} above. Then,
applying Pisier's characterization \cite{P0} of
$\mathrm{K}$-convexity, conditions i) and iii) are equivalent. Now
we prove the equivalence between iii) and iv). To that aim we can
fix $1 < p < \infty$ without lost of generality (note that iii) is
independent of the index $p \in (1, \infty)$ by its equivalence
with i) and Remark \ref{Remark-p-Independence}). The implication
iii) $\Rightarrow$ iv) is trivial. Reciprocally, let us assume
that $S_p(E)$ contains $\ell_1^n$'s uniformly. Note that, since
$\ell_1^n$ carries the maximal operator space structure, any
Banach space embedding of $\ell_1^n$ is automatically a
semi-complete embedding with the same constants. Then, Corollary
\ref{Corollary-Semi-Embedding} claims that the family $$\Big\{
S_1^n \, \Big| \ n \ge 1 \Big\}$$ also embeds semi-completely
uniformly in $S_p(E)$. That is, there exists $\mathrm{c}
> 1$ and embeddings $\Lambda_n: S_1^n \rightarrow S_p(E)$ such
that $$\|\Lambda_n\|_{cb} \|\Lambda_n^{-1}\| \le \mathrm{c}.$$
Now, given a finite subset $\Gamma$ of $\Sig$, we also consider
the map $$\mathrm{S}_{\Gamma}: A \in \Lebesgue_1(\Gamma)
\longmapsto \bigoplus_{\sig \in \Gamma} \Degree A^{\sig} \in
S_1^{\mathrm{N}} \qquad \mbox{for} \qquad \mathrm{N} = \sum_{\sig
\in \Gamma} \Degree.$$ Finally, let $\mathrm{R}_{\Gamma}:
\Lebesgue_1(\Gamma) \to S_p(E)$ stand for $\Lambda_{\mathrm{N}}
\circ \mathrm{S}_{\Gamma}$. Then we have
$$\|\mathrm{R}_{\Gamma}\|_{cb} \|\mathrm{R}_{\Gamma}^{-1}\| \le
\|\Lambda_{\mathrm{N}}\|_{cb} \|\Lambda_{\mathrm{N}}^{-1}\| \le
\mathrm{c},$$ since $\mathrm{S}_{\Gamma}$ is a complete isometry.
In summary, the $\Lebesgue_1(\Gamma)$'s embed semi-completely
uniformly in $S_p(E)$. This proves the implication iv)
$\Rightarrow$ iii). It remains to see that ii) is equivalent to
some (any) of the other conditions. As in the commutative case,
the implication ii) $\Rightarrow$ iv) follows from Remark
\ref{Remark-Uniform-lambda} and by plugging in the \lq${}$right
unit vectors\rq, for details see \cite{Pa}. The converse
implication iv) $\Rightarrow$ ii) (a bit more technical) is the
main result in \cite{Pa}. This completes the proof. \fin

\begin{remark}
\emph{Theorem \ref{Theorem-Pisier} implies the $\Sig$-independence
of $\mathrm{OB}_{\Sig}$-convexity.}
\end{remark}

\begin{remark}
\emph{We have already mentioned that semi-complete and Banach
space embeddings of $\ell_1^n$'s are the same since $\ell_1^n$
carries the maximal o.s.s. It is worthy of mention that, although
$\Lebesgue_1(\Gamma)$'s are not longer equipped with the maximal
operator space structure, a similar property holds for the latter
spaces. Indeed, it is clear that if $\Lebesgue_1(\Gamma)$'s are
uniformly contained in $S_p(E)$ in the Banach space sense, then
$S_p(E)$ also contains $\ell_1^n$'s uniformly. Finally, by Theorem
\ref{Theorem-Pisier} we see that $\Lebesgue_1(\Gamma)$'s embed
semicompletely uniformly in $S_p(E)$. The converse is trivial.}
\end{remark}

\section{Operator space type and cotype}
\label{Section6}

The notions of Fourier type and cotype of an operator space with
respect to a noncommutative compact group were already defined in
the Introduction. These are particular cases of a more general
notion of type and cotype for operator spaces introduced in
\cite{GP2}. In that paper, the (\emph{uniformly bounded})
\emph{quantized orthonormal systems} play the same role of the
uniformly bounded orthonormal systems in the classical theory.
Some relevant examples of this notion are the dual object of a
noncommutative compact group and the quantized analog of the
Steinhaus system introduced in \cite{MP}. Before introducing the
notions of type and cotype for operator spaces, let us recover the
classical notions. Let $\varepsilon_1, \varepsilon_2, \ldots$ be a
sequence of random signs or independent Steinhaus variables over a
probability space $(\Omega, \mathsf{A}, \mu)$. Given $1 \le p \le
2$, a Banach space $\mathrm{X}$ is called of type $p$ when there
exists a constant $\mathsf{T}_p(\mathrm{X})$ such that $$\Big(
\int_{\Omega} \Big\| \sum_{k=1}^n x_k \varepsilon_k(\omega)
\Big\|_{\mathrm{X}}^{p'} d\mu(\omega) \Big)^{1/p'} \le
\mathsf{T}_p(\mathrm{X}) \, \Big( \sum_{k=1}^n
\|x_k\|_{\mathrm{X}}^p \Big)^{1/p}$$ for any finite family $x_1,
x_2, \ldots, x_n$ in $\mathrm{X}$. As we mentioned in the
Introduction, the basic idea is to replace the random variables
$(\varepsilon_k)$ by a sequence $U_1, U_2, \ldots$ of independent
random unitaries. That is, each $U_k: \Omega \to U(d_k)$ is a
random unitary $d_k \times d_k$ matrix uniformly distributed in
the unitary group $U(d_k)$ with respect to the normalized Haar
measure. In this setting, we might define the following notion of
type $$\Big( \int_{\Omega} \Big\| \sum_{k=1}^n d_k \sum_{i,j=1}^n
A_k(i,j) U_k(j,i) \Big\|_{\mathrm{X}}^{p'} d\mu(\omega)
\Big)^{1/p'} \le \widetilde{\mathsf{T}}_p(\mathrm{X}) \, \Big(
\sum_{k=1}^n d_k \|A_k\|_{S_p^{d_k}(\mathrm{X})}^p \Big)^{1/p}.$$
We want to point out that the right hand side is only well-defined
for operator spaces. Moreover, this notion depends on the
dimension $d_k$ and their multiplicity. Note that the presence of
$d_k$'s in the inequality stated above is quite natural in view of
the Peter-Weyl theorem and the connection (explained in the
Introduction) with the Hausdorff-Young inequality for non-abelian
compact groups. Let us give the precise definitions. The
\emph{quantized Steinhaus system} associated to $(\Sig,
\mathbf{d}_{\Sig})$ is defined as a collection $$\mathbf{S}_{\Sig}
= \Big\{ \zeta^{\sig}: \Omega \rightarrow U(\Degree) \Big\}_{\sig
\in \Sig}$$ of independent uniformly distributed random unitaries
with respect to the set of parameters $(\Sig, \mathbf{d}_{\Sig})$.
Given an operator space $E$ and a function $f \in L_2(\Omega;E)$,
we can consider the Fourier coefficients of $f$ with respect to
the quantized Steinhaus system $$\widehat{f}_{\mathbf{S}}(\sig) =
\int_{\Omega} f(\omega) \zeta^{\sig}(\omega)^{\ast} d
\mu(\omega).$$ Let $\mathbf{St}_p(\Sig;E)$ be the closure in
$L_p(\Omega;E)$ of the subspace given by functions $$f_{\Gamma} =
\sum_{\sig \in \Gamma} \Degree \mbox{tr} (A^{\sig} \zeta^{\sig})
\qquad \mbox{with} \qquad A^{\sig} \in M_{\Degree} \otimes E$$ and
$\Gamma$ a finite subset of $\Sig$. We shall write
$\mathbf{St}_p(\Sig)$ for the scalar-valued case. Then, given $1
\le p \le 2$, we say that the operator space $E$ has
\textbf{$\Sig$-type} $p$ when the following inequality holds for
any function $f \in \mathbf{St}_{p'}(\Sig;E)$ $$\Big(
\int_{\Omega} \|f(\omega)\|_E^{p'} \ d\mu(\omega) \Big)^{1/p'}
\le_{cb} \mathcal{K}_p^1(E,\mathbf{S}_{\Sig}) \, \Big( \sum_{\sig
\in \Sig} \Degree
\|\widehat{f}_{\mathbf{S}}(\sig)\|_{S_p^{\Degree}(E)}^p
\Big)^{1/p}.$$ In a similar way, \textbf{$\Sig$-cotype} $p'$ means
that any $f \in \mathbf{St}_p(\Sig;E)$ satisfies $$\Big(
\sum_{\sig \in \Sig} \Degree
\|\widehat{f}_{\mathbf{S}}(\sig)\|_{S_{p'}^{\Degree}(E)}^{p'}
\Big)^{1/p'} \le_{cb} \mathcal{K}_{p'}^2(E,\mathbf{S}_{\Sig}) \,
\Big( \int_{\Omega} \|f(\omega)\|_E^p \ d\mu(\omega)
\Big)^{1/p}.$$ Recall that the symbol $\le_{cb}$ means the
complete boundedness of the corresponding linear map. The best
constants $\mathcal{K}_p^1(E,\mathbf{S}_{\Sig})$ and
$\mathcal{K}_{p'}^2(E,\mathbf{S}_{\Sig})$ in the inequalities
stated above are called the $\Sig$-type $p$ and $\Sig$-cotype $p'$
constants of $E$. More concretely, using the spaces
$\Lebesgue_p(\Sig;E)$ introduced in Section \ref{Section5}, the
given definitions of $\Sig$-type and $\Sig$-cotype can be
rephrased be requiring the complete boundedness of the following
operators $$\mathrm{T}_p: A \in \Lebesgue_p(\Sig;E) \longmapsto
\sum_{\sig \in \Sig} \Degree \mbox{tr}(A^{\sig} \zeta^{\sig}) \in
\mathbf{St}_{p'}(\Sig;E),$$ $$\mathrm{C}_{p'}: \sum_{\sig \in
\Sig} \Degree \mbox{tr}(A^{\sig} \zeta^{\sig}) \in
\mathbf{St}_p(\Sig;E) \longmapsto A \in \Lebesgue_{p'}(\Sig;E).$$

\begin{remark} \label{Remark-Khintchine-Kahane}
\emph{Let us recall that $\Sig_0$ stands for the commutative set
of parameters defined in Section \ref{Section5}. The classical
Khintchine inequalities can be rephrased by saying that the norm
of $\mathbf{St}_p(\Sig_0)$, regarded as a Banach space, is
equivalent to that of $\mathbf{St}_q(\Sig_0)$ whenever $1 \le p
\neq q < \infty$. On the other hand, by means of the
noncommutative Khintchine inequalities \cite{Lu,LuP}, it turns out
that the norm of $\mathbf{St}_p(\Sig_0)$ is \emph{not} completely
equivalent to that of $\mathbf{St}_q(\Sig_0)$. That is, the
operator spaces $\mathbf{St}_p(\Sig_0)$ and
$\mathbf{St}_q(\Sig_0)$ are isomorphic but not completely
isomorphic. More generally, $\mathbf{St}_p(\Sig)$ is Banach
isomorphic but not completely isomorphic to $\mathbf{St}_q(\Sig)$,
see \cite{MP} for the details. Therefore, each space
$\mathbf{St}_q(\Sig)$ in the definition of $\Sig$-type $p$ and
$\Sig$-cotype $p'$ gives a priori a different notion!}
\end{remark}

\begin{remark} \label{Remark-Non-Trivial}
\emph{As in the classical theory, every operator space has
$\Sigma$-type $1$ and $\Sigma$-cotype $\infty$. An operator space
$E$ has \emph{non-trivial $\Sig$-type} whenever it has $\Sig$-type
$p$ for some $1 < p \le 2$. According to \cite{Pa} and in contrast
with the commutative theory, $\mathrm{OK}$-convexity is not
equivalent to having non-trivial $\Sig$-type. Indeed, the operator
Hilbert spaces $R$ and $C$ fail this equivalence since both are
$\mathrm{OK}$-convex operator spaces but do not have $\Sig$-type
for any $1 < p \le 2$. This constitutes an important difference
between the classical and the noncommutative contexts. Namely, it
turns out that we can not expect an operator space version of the
Maurey-Pisier theorem \cite{MauP} since the simplest form of this
result asserts that the property of having non-trivial type is
equivalent to $\mathrm{K}$-convexity.}
\end{remark}

\begin{remark}
\emph{It is not clear whether or not the notions of $\Sig$-type
and $\Sig$-cotype depend on $(\Sig, \mathbf{d}_{\Sig})$. Moreover,
if we replace the quantized Steinhaus system by the dual object of
a noncommutative compact group $\mathrm{G}$, we can ask ourselves
the same question for the notions of Fourier type and cotype. Note
that this \emph{group independence} is an open problem even in the
commutative theory. The reader is referred to the paper \cite{HL}
for more information on this problem.}
\end{remark}

The $\Sig$-type (resp. $\Sig$-cotype) becomes a stronger condition
on any operator space as the exponent $p$ (resp. $p'$) approaches
$2$. In particular, given an operator space $E$ we consider (as in
the Banach space context) the notions of \emph{sharp} $\Sig$-type
of $E$ (i.e. the supremum over all $1 \le p \le 2$ for which $E$
has $\Sig$-type $p$) as well as \emph{sharp} $\Sig$-cotype of $E$
(i.e. the infimum over all $2 \le p' \le \infty$ for which $E$ has
$\Sig$-type $p'$). The aim of this section is to investigate the
sharp $\Sig$-type and $\Sig$-cotype indices of Lebesgue spaces,
either commutative or not. However, as we shall see below, some
other related problems will be solved with the same techniques.

\subsection{Sharp $\Sig$-type of $L_p$ for $1 \le p \le 2$}

We begin with the finite dimensional $\Sig$-type constants for any
\emph{bounded} set of parameters $(\Sig, \mathbf{d}_{\Sig})$. More
concretely, let us consider a set of parameters $(\Sig,
\mathbf{d}_{\Sig})$ with $\mathbf{d}_{\Sig}$ bounded. Then, given
a finite subset $\Gamma$ of $\Sig$, we shall write
$\ell_p(\Gamma)$ to denote the space of functions $\xi: \Gamma \to
\C$ endowed with the customary norm $$\|\xi\|_{\ell_p(\Gamma)} =
\Big( \sum_{\sig \in \Gamma} |\xi(\sig)|^p \Big)^{1/p}.$$ Let us
consider the function $f: \Omega \to \ell_p(\Gamma)$ defined by
$$f = \sum_{\sig \in \Gamma} \Degree \mbox{tr}
(\widehat{f}_{\mathbf{S}}(\sig) \zeta^{\sig}) \qquad \mbox{with}
\qquad \widehat{f}_{\mathbf{S}}(\sig) = e_{11} \otimes
\delta_{\sigma} \in M_{\Degree} \otimes \ell_p(\Gamma).$$ Then we
recall that $$\Big( \int_{\Omega} |\sqrt{\Degree}
\zeta_{11}^{\sig} |^q d\mu \Big)^{\frac1q} \thicksim \Big(
\int_{\Omega} |\sqrt{\Degree} \zeta_{11}^{\sig} |^2 d\mu
\Big)^{\frac12} = 1, \qquad \mbox{for any} \quad 1 \le q <
\infty.$$ Indeed, the norm equivalence follows from the analog of
the Khintchine-Kahane inequalities for the quantized Steinhaus
system, proved in \cite{MP}. The last equality follows from the
definition of $\mathbf{S}_{\Sig}$. In particular, if we use the
symbol $\lesssim$ to denote an inequality up to a universal
positive constant, then we have the following estimate for any $1
\le p < q \le 2$
\begin{eqnarray*}
|\Gamma|^{1/p} & \lesssim & \Big( \int_{\Omega} \sum_{\sig \in
\Gamma} |\Degree \zeta_{11}^{\sig}|^p d \mu \Big)^{1/p} \\ & \le &
\Big( \int_{\Omega} \Big\| \sum_{\sig \in \Gamma} \Degree
\mbox{tr} (\widehat{f}_{\mathbf{S}}(\sig) \zeta^{\sig})
\Big\|_{\ell_p(\Gamma)}^{q'} d \mu \Big)^{1/q'} \\ & \le &
\mathcal{K}_q^1(\ell_p(\Gamma), \mathbf{S}_{\Sig}) \, \Big(
\sum_{\sig \in \Gamma} \Degree
\|\widehat{f}_{\mathbf{S}}(\sig)\|_{S_q^{\Degree}(\ell_p(\Gamma))}^q
\Big)^{1/q} \lesssim \ \mathcal{K}_q^1(\ell_p(\Gamma),
\mathbf{S}_{\Sig}) \, |\Gamma|^{1/q}.
\end{eqnarray*}
In other words $$c \, |\Gamma|^{1/p - 1/q} \le
\mathcal{K}_q^1(\ell_p(\Gamma), \mathbf{S}_{\Sig}) \le
|\Gamma|^{1/p - 1/q},$$ for some constant $0 < c \le 1$. The upper
estimate is much simpler and it can be found in \cite{GP1}.
Therefore, since any infinite dimensional (either commutative or
noncommutative) $L_p$ space contains completely isometric copies
of $\ell_p(\Gamma)$ for any finite subset $\Gamma$ of $\Sig$, we
deduce that any infinite dimensional $L_p$ space has sharp
$\Sig$-type $p$ for any bounded set of parameters $(\Sig,
\mathbf{d}_{\Sig})$. However, it is evident that our argument
doesn't work for unbounded sets of parameters. This case requires
to find the right matrices which give the optimal constants. In
the following theorem we compute the finite dimensional constants
for the Schatten classes.

\begin{theorem} \label{Theorem-Sigma-Type}
If $1 \le p < q \le 2$, the estimate
$$\mathcal{K}_q^1(S_p^{d_{\sig}}, \mathbf{S}_{\Sig}) \ge
d_{\sig}^{2(1/p - 1/q)}$$ holds for any unbounded set of
parameters $(\Sig, \mathbf{d}_{\Sig})$ and any element $\sig$ of
$\Sig$.
\end{theorem}

\dem Let us take $f: \Omega \to S_p^{\Degree}$ so that
$\widehat{f}_{\mathbf{S}}(\xi) = 0$ if $\xi \in \Sig \setminus
\{\sig\}$ and $$\widehat{f}_{\mathbf{S}}(\sig) = \Big(
\sum_{i=1}^{d_{\sig}} e_{i1} \otimes e_{i1} \Big) \otimes \Big(
\sum_{j=1}^{d_{\sig}} e_{1j} \otimes e_{1j} \Big) \in
C_q^{d_{\sig}} \otimes_h C_p^{d_{\sig}} \otimes_h R_p^{d_{\sig}}
\otimes_h R_q^{d_{\sig}} = S_q^{d_{\sig}} (S_p^{d_{\sig}}).$$
Then, the following estimate holds by definition of $\Sig$-type
$$\Big( \int_{\Omega} \Big\| d_{\sig} \mbox{tr}
(\widehat{f}_{\mathbf{S}}(\sig) \zeta^{\sig})
\Big\|_{S_p^{d_{\sig}}}^{q'} d\mu \Big)^{1/q'} \le
\mathcal{K}_q^1(S_p^{d_{\sig}},\mathbf{S}_{\Sig}) \,
d_{\sig}^{1/q}
\|\widehat{f}_{\mathbf{S}}(\sig)\|_{S_q^{d_{\sig}}(S_p^{d_{\sig}})}.$$
Note that we have $$\mbox{tr}(\widehat{f}_{\mathbf{S}}(\sig)
\zeta^{\sig}) = \sum_{i,j = 1}^{d_{\sig}} e_{ij} \otimes
\zeta_{ji}^{\sig} = (\zeta^{\sig})^{\mbox{t}}.$$ Thus, since the
$\zeta^{\sig}$'s are unitary, the left hand side of the inequality
above is $\Degree^{1 + \frac1p}$. On the other hand, we need to
compute the norm of $\widehat{f}_{\mathbf{S}}(\sig)$ in
$S_q^{\Degree}(S_p^{\Degree})$. Since the Haagerup tensor product
commutes with complex interpolation, it is not difficult to check
that the following natural identifications are Banach space
isometries $$C_q^{d_{\sig}} \otimes_h C_p^{d_{\sig}} =
S_r^{d_{\sig}} = R_p^{d_{\sig}} \otimes_h R_q^{d_{\sig}} \qquad
\mbox{with} \qquad \frac{1}{r} = \frac12 \Big( 1 - \frac1p +
\frac1q \Big).$$ For instance, $$C_{\infty}^{\Degree} \otimes_h
C_p^{\Degree} = [C_{\infty}^{\Degree} \otimes_h
C_{\infty}^{\Degree}, C_{\infty}^{\Degree} \otimes_h
C_1^{\Degree}]_{1/p} = [S_2^{\Degree}, S_{\infty}^{\Degree}]_{1/p}
= S_{2p'}^{\Degree},$$ $$C_q^{d_{\sig}} \otimes_h C_p^{d_{\sig}}=
[C_{\infty}^{\Degree} \otimes_h C_p^{\Degree}, C_p^{\Degree}
\otimes_h C_p^{\Degree}]_{p/q} =
[S_{2p'}^{\Degree},S_2^{\Degree}]_{p/q} = S_r^{\Degree}.$$ In
particular, due to our choice of $\widehat{f}_{\mathbf{S}}(\sig)$,
we can write
$$\|\widehat{f}_{\mathbf{S}}(\sig)\|_{S_q^{\Degree}(S_p^{\Degree})}
= \|1_{M_{\Degree}}\|_{S_r^{\Degree}}^2 = \Degree^{2/r}.$$
Combining our previous results, we obtain the desired estimate.
\fin

\begin{corollary} \label{Corollary-Sigma-Type}
If $1 \le p < q \le 2$, the estimate
$$\mathcal{K}_q^1(\ell_p^{d_{\sig}^2}, \mathbf{S}_{\Sig}) \gtrsim
d_{\sig}^{2(1/p - 1/q)}$$ holds for any unbounded set of
parameters $(\Sig, \mathbf{d}_{\Sig})$ and any element $\sig$ of
$\Sig$.
\end{corollary}

\dem By Theorem \ref{Theorem-Embedding}, we have
$$\mathcal{K}_q^1(S_p^{d_{\sig}}, \mathbf{S}_{\Sig}) \lesssim \,
\mathcal{K}_q^1(S_q(\ell_p^{d_{\sig}^2}), \mathbf{S}_{\Sig}) \le
\mathcal{K}_q^1(\ell_p^{d_{\sig}^2}, \mathbf{S}_{\Sig}).$$ The
last inequality follows by Minkowski inequality for operator
spaces, see \cite{GP1}. \fin

\begin{remark} \label{Remark-Type-Cotype}
\emph{The arguments applied up to now also provide the finite
dimensional estimates for the $\Sig$-cotype constants when $2 \le
q' < p' \le \infty$. Namely, the following estimates hold
$$\mathcal{K}_{q'}^2(S_{p'}^{\Degree}, \mathbf{S}_{\Sig}) \ge
\Degree^{2(1/q' - 1/p')} \qquad \mbox{and} \qquad
\mathcal{K}_{q'}^2(\ell_{p'}^{\Degree^2}, \mathbf{S}_{\Sig})
\gtrsim \Degree^{2(1/q' - 1/p')}.$$}
\end{remark}

\begin{remark}
\emph{By a simple result of \cite{GP1}, we have
$\mathcal{K}_q^1(S_p^{\Degree}, \mathbf{S}_{\Sig}) \le d_{cb}
(S_p^{\Degree}, S_q^{\Degree})$. In particular, in Theorem
\ref{Theorem-Sigma-Type} we actually have equality
$$\mathcal{K}_q^1(S_p^{\Degree}, \mathbf{S}_{\Sig}) =
\Degree^{2(1/p - 1/q)}.$$ A similar argument applies to Corollary
\ref{Corollary-Sigma-Type}. In summary, our estimates provide the
\emph{exact} order of growth of the $\Sig$-type (resp.
$\Sig$-cotype by Remark \ref{Remark-Type-Cotype}) constants of the
corresponding finite dimensional Lebesgue spaces considered above.
Moreover, now we can prove the claim given in Remark
\ref{Remark-Optimal}. Namely, let us consider the set $\Sig = \N$
with $d_k = k$ for all $k \ge 1$. Then, if $\Psi_{pq}: S_q^n
\rightarrow S_p(\ell_q^m)$ is a $cb$ embedding with constants not
depending on the dimensions $n$ and $m$, Corollary
\ref{Corollary-Sigma-Type} provides the following estimate}
\begin{eqnarray*}
\mathcal{K}_q^1(S_p^n, \mathbf{S}_{\Sig}) & \le &
\|\Psi_{pq}\|_{cb} \|\Psi_{pq}^{-1}\|_{cb} \,
\mathcal{K}_q^1(S_q(\ell_p^m), \mathbf{S}_{\Sig}) \\ & \le &
\|\Psi_{pq}\|_{cb} \|\Psi_{pq}^{-1}\|_{cb} \,
\mathcal{K}_q^1(\ell_p^m, \mathbf{S}_{\Sig}) \\ & \le &
\|\Psi_{pq}\|_{cb} \|\Psi_{pq}^{-1}\|_{cb} \, m^{1/p - 1/q}.
\end{eqnarray*}
\emph{Therefore, since $\mathcal{K}_q^1(S_p^n, \mathbf{S}_{\Sig})
= n^{2(1/p - 1/q)}$, we conclude by taking $n$ arbitrary large.}
\end{remark}

\begin{remark}
\emph{The main topic of \cite{GMP} is the sharp Fourier type and
cotype of $L_p$ spaces. Given $1 \le p \le 2$, it is showed that
$L_p$ has sharp Fourier type $p$ with respect to any compact
semisimple Lie group. The arguments employed are very different.
Namely, the key point is a Hausdorff-Young type inequality for
functions defined on a compact semisimple Lie group with arbitrary
small support. However, the sharp Fourier cotype of $L_p$ for $1
\le p \le 2$ is left open in \cite{GMP}. Now we can solve it by
using Corollary \ref{Corollary-Sigma-Type} and the following
inequality $$\mathcal{K}_{q'}^2(L_p, \widehat{\mathrm{G}}) \ge
\mathcal{K}_q^1(L_p, \mathbf{S}_{\widehat{\mathrm{G}}}).$$ Here
$\mathcal{K}_{q'}^2(L_p, \widehat{\mathrm{G}})$ denotes the
Fourier cotype $q'$ constant of $L_p$ with respect to $\mathrm{G}$
and $\mathbf{S}_{\widehat{\mathrm{G}}}$ stands for the quantized
Steinhaus system with the parameters given by the degrees of the
irreducible representations of $\mathrm{G}$. That inequality is a
particular case of the noncommutative version of the contraction
principle given in \cite{MP}. This solves the problem posed in
\cite{GMP} not only for compact semisimple Lie groups, but for any
non-finite topological compact group.}
\end{remark}

\subsection{Sharp $\Sig$-cotype of $L_p$ for $1 \le p \le 2$}

Given any $\sigma$-finite measure space $(\widetilde{\Omega},
\mathsf{B}, \nu)$, any set of parameters $(\Sig,
\mathbf{d}_{\Sig})$ and any finite subset $\Gamma$ of $\Sig$, let
us consider a family of matrices $$\mathbf{A} = \Big\{A^{\sig} \in
M_{\Degree} \otimes L_p(\widetilde{\Omega}) \Big\}_{\sig \in
\Gamma}.$$ Then, we can estimate the norm of $\mathbf{A}$ in
$\Lebesgue_2(\Sig;L_p(\widetilde{\Omega}))$ for any $1 \le p \le
2$ as follows. First, Minkowski inequality and Plancherel theorem
give
\begin{eqnarray*}
\Big( \sum_{\sig \in \Gamma} \Degree
\|A^{\sig}\|_{S_2^{\Degree}(L_p(\widetilde{\Omega}))}^2
\Big)^{1/2} & \le & \Big( \int_{\widetilde{\Omega}} \Big[
\sum_{\sig \in \Gamma} \Degree \|A^{\sig}(x)\|_{S_2^{\Degree}}^2
\Big]^{p/2} d \nu(x) \Big)^{1/p}
\\ & = & \Big( \int_{\widetilde{\Omega}} \Big[ \int_{\Omega} \Big| \sum_{\sig
\in \Gamma} \Degree \mbox{tr} (A^{\sig} \zeta^{\sig}) \Big|^2 d
\mu \Big]^{p/2} d \nu \Big)^{1/p}
\end{eqnarray*}
Second, by the analog given in \cite{MP} of Khintchine-Kahane
inequalities for $\mathbf{S}_{\Sig}$ $$\Big(
\int_{\widetilde{\Omega}} \Big[ \int_{\Omega} \Big| \sum_{\sig \in
\Gamma} \Degree \mbox{tr} (A^{\sig} \zeta^{\sig}) \Big|^2 d \mu
\Big]^{\frac{p}{2}} d \nu \Big)^{\frac1p} \thicksim \Big(
\int_{\Omega} \Big[ \int_{\widetilde{\Omega}} \Big| \sum_{\sig \in
\Gamma} \Degree \mbox{tr} (A^{\sig} \zeta^{\sig}) \Big|^p d \nu
\Big]^{\frac2p} d \mu \Big)^{\frac12}.$$ Therefore, there exists
some constant $\mathrm{c}$ such that $$\Big( \sum_{\sig \in
\Gamma} \Degree
\|A^{\sig}\|_{S_2^{\Degree}(L_p(\widetilde{\Omega}))}^2
\Big)^{1/2} \le \mathrm{c} \, \Big( \int_{\Omega} \Big\|
\sum_{\sig \in \Gamma} \Degree \mbox{tr} (A^{\sig} \zeta^{\sig}
(\omega)) \Big\|_{L_p(\widetilde{\Omega})}^2 d \mu (\omega)
\Big)^{1/2}$$ for any family of matrices $\mathbf{A}$. In other
words, we have proved that the mapping $\mathrm{C}_2$ defined
above is bounded when we take values in $L_p(\widetilde{\Omega})$.
However, we can not claim $\Sig$-cotype $2$ unless we prove that
the same operator $\mathrm{C}_2$ is not only bounded, but
completely bounded. Now, looking at Remark
\ref{Remark-Khintchine-Kahane}, we realize that our arguments do
not work to show the complete boundedness. In this paragraph we
study this problem. We begin by computing the sharp cotype of
$S_q(S_p)$ as a Banach space. This will be the key to find the
sharp $\Sig$-cotype indices of $L_p$ spaces. We want to point out
that this fact was independently discovered by Lee in \cite{Le}.

\begin{lemma} \label{Lemma-SqSp}
The Schatten class $S_q(S_p)$ has sharp Banach cotype $r$ with
$$\frac1r = \frac12 \Big( 1 - \frac1p + \frac1q \Big) \qquad
\mbox{whenever} \qquad 1 \le p \le 2  \quad \mbox{and} \quad p \le
q \le p'.$$
\end{lemma}

\dem First, we show that $S_q(S_p)$ has cotype $r$. The case $p
> 1$ is simple. Indeed, we just need to check that the predual
$S_{q'}(S_{p'})$ has Banach type $r'$. To that aim we observe that
$$S_{q'}(S_{p'}) = [S_p(S_{p'}),S_{p'}(S_{p'})]_\theta \qquad
\mbox{with} \qquad 1 - \frac1q = \frac{1-\theta}{p} + \theta
\Big(1 - \frac1p \Big).$$ Moreover, we have $$S_p(S_{p'}) =
[S_2(S_2),S_1(S_{\infty})]_{\eta} \qquad \mbox{with} \qquad
\frac{1}{p} = \frac{1-\eta}{2} + \frac{\eta}{1}.$$ Hence
$S_p(S_{p'})$ has type $p$ and, since $S_{p'}(S_{p'})$ has type
$2$, $S_{q'}(S_{p'})$ has type $s$ with $$\frac{1}{s} =
\frac{1-\theta}{p} + \frac{\theta}{2} = 1 - \frac1q + \theta
\Big(\frac1p - \frac12 \Big) = 1 - \frac1q + \frac12 \Big( \frac1p
-1 + \frac1q \Big) = 1 - \frac1r.$$ It remains to see that
$S_q(S_1)$ has cotype $2q$. Let us denote by $\mathcal{R}_p$ the
subspace generated in $L_p(\Omega)$ by the sequence $\rad_1,
\rad_2, \ldots$ of Rademacher functions. Then, if
$\mathcal{R}_p(E)$ stands for the closure of the tensor product
$\mathcal{R}_p \otimes E$ in $L_p(\Omega;E)$, we need to see that
the following mapping is bounded $$\mathrm{C}_{2q}: \sum_{k=1}^n
\rad_k \otimes x_k \in \mathcal{R}_2(S_q(S_1)) \longmapsto
\sum_{k=1}^n \delta_k \otimes x_k \in \ell_{2q}(S_q(S_1)).$$ First
we recall that, according to Khintchine-Kahane and Minkowski
inequalities, the following natural map is contractive
$$\mathcal{R}_2(S_q(S_1)) \simeq S_q(\mathcal{R}_q(S_1))
\rightarrow S_q(\mathcal{R}_1(S_1)).$$ By the well-known complete
isomorphism $\mathcal{R}_1 \simeq R+C$, which follows from the
noncommutative Khintchine inequalities (see \cite{LuP,P2}), we can
write $S_q(\mathcal{R}_1(S_1))$ as the sum $S_q(S_1(R)) +
S_q(S_1(C))$. Therefore, it suffices to see that the following
natural mappings $$\begin{array}{lrcl} \mathrm{S}: & S_q(S_1(R)) &
\rightarrow & \ell_{2q}(S_q(S_1)) \\ \mathrm{T}: & S_q(S_1(C)) &
\rightarrow & \ell_{2q}(S_q(S_1)), \end{array}$$ which send the
canonical basis of $R$ or $C$ to the canonical basis of
$\ell_{2q}$, are bounded. Since both cases are similar, we only
prove the boundedness of $\mathrm{T}$. To that aim we recall that,
since $S_q(S_1(C)) = [S_{\infty}(S_1(C)), S_1(S_1(C))]_{1/q}$, it
suffices to prove the boundedness of $$\begin{array}{lrcl}
\mathrm{T}_0: & S_{\infty}(S_1(C)) & \rightarrow &
\ell_{\infty}(S_{\infty}(S_1))
\\ \mathrm{T}_1: & S_1(S_1(C)) & \rightarrow & \ell_2(S_1(S_1)).
\end{array}$$ If we observe that $\mathrm{T}_0$ factors through
$S_{\infty}(S_1(\ell_{\infty}))$, it is clear that $\mathrm{T}_0$
is even contractive. To show that $\mathrm{T}_1$ is bounded, let
us consider a finite family $x_1,x_2,.., x_n$ of elements in
$S_1(S_1)$. Then, since $S_1(S_1(C))$ embeds completely
isometrically in $S_1(\N^3)$, we know from \cite{TJ} that it has
Banach cotype $2$ so that we get
\begin{eqnarray*}
\Big\| \sum_{k=1}^n \delta_k \otimes x_k \Big\|_{\ell_2(S_1(S_1))}
& = & \Big( \sum_{k=1}^n \|x_k \otimes e_{k1}\|_{S_1(S_1(C))}^2
\Big)^{1/2} \\ & \le & \mathrm{c} \, \int_0^1 \Big\| \sum_{k=1}^n
\rad_k(t) (x_k \otimes e_{k1}) \Big\|_{S_1(S_1(C))} dt \\ & = &
\mathrm{c} \ \Big\| \sum_{k=1}^n x_k \otimes e_{k1}
\Big\|_{S_1(S_1(C))}.
\end{eqnarray*}
The last equality follows since $$\Big\| \sum_{k=1}^n \rad_k(t)
(x_k \otimes e_{k1}) \Big\|_{S_1(S_1(C))} = \Big\| \Big(
\sum_{k=1}^n (\rad_k(t) x_k)^* (\rad_k(t) x_k) \Big)^{1/2}
\Big\|_{S_1(S_1)}.$$ This gives the boundedness of $\mathrm{T}_1$
and consequently the map $\mathrm{C}_{2q}$ is also bounded. In
summary, we have seen that $S_q(S_p)$ has Banach cotype $r$ in the
range of parameters considered. To complete the proof, we need to
see that this exponent is sharp. However, recalling that
$$S_q(S_p) = C_q\otimes_h C_p \otimes_h R_p \otimes_h R_q,$$ we
can regard $C_q \otimes_h C_p$ as a subspace of $S_q(S_p)$. Now,
since the Haagerup tensor product commutes with complex
interpolation, we obtain the following Banach space isometries
$$C_q \otimes_h C_p = [C_{p'} \otimes_h C_p, C_p \otimes_h
C_p]_{\theta} = [S_{p'},S_2]_{\theta} \qquad \mbox{with} \qquad
\frac1q = 1 - \frac1p + \theta \Big( \frac2p - 1 \Big).$$ This
gives that $C_q \otimes_h C_p = S_r$ as a Banach space, we leave
the details to the reader. Therefore, $S_q(S_p)$ can not have
better cotype than $r$. This completes the proof. \fin

Let us recall that the commutative set of parameters $(\Sig_0,
\mathbf{d}_{\Sig_0})$ is the given by $\Sig_0 = \N$ where we take
$\Degree = 1$ for all $\sig \in \Sig_0$. In the following result
we show that, in contrast with the Banach space situation, any
infinite dimensional (commutative or noncommutative) $L_p$ space
with $p \neq 2$ fails to have $\Sig$-cotype $2$.

\begin{theorem} \label{Theorem-Sharp-Cotype}
Any infinite dimensional $L_p$ space has sharp $\Sig$-cotype
$\max(p,p')$.
\end{theorem}

\demone As it was pointed out in \cite{GP2}, it is obvious the any
$L_p$ space has $\Sig$-cotype $\max(p,p')$ with respect to any set
of parameters $\Sig$. Let us see that this exponent is sharp when
$\mathbf{d}_{\Sig}$ is bounded. In this particular case, it
clearly suffices to consider the commutative set of parameters
$\Sig_0$. We also assume that $1 \le p \le 2$ since the case $2
\le p \le \infty$ has been considered in Remark
\ref{Remark-Type-Cotype}. Moreover, since any infinite dimensional
$L_p$ space contains a completely isometric copy of $\ell_p$, it
suffices to check it for $\ell_p$. Now, let us assume that
$\ell_p$ has $\Sig_0$-cotype $q'$, for some $q' < p'$. Then we can
argue as in Corollary \ref{Corollary-Sigma-Type}. Namely,
combining Theorem \ref{Theorem-Embedding} with Minkowski
inequality for operator spaces, we have
$$\mathcal{K}_{q'}^2(S_{q'}(S_p); \mathbf{S}_{\Sig_0}) \lesssim
\mathcal{K}_{q'}^2(S_{q'}(\N^2;\ell_p); \mathbf{S}_{\Sig_0}) \le
\mathcal{K}_{q'}^2(\ell_p; \mathbf{S}_{\Sig_0}).$$ Now, by Lemma
\ref{Lemma-SqSp}, the best $\Sig_0$-cotype we can expect to have
is $r$ where $$\frac1r = \frac{1}{2p'} + \frac{1}{2q'} <
\frac{1}{q'}.$$ Therefore, we deduce that $r > q'$ and the result
follows by contradiction. \fin

\demtwo Arguing as in the previous case, it suffices to see that
$\ell_p$ has sharp $\Sig$-cotype $p'$ for $1 \le p \le 2$. Let us
assume that $\ell_p$ has $\Sig$-cotype $q'$ for some $q' < p'$.
Then, again by Theorem \ref{Theorem-Embedding} and Minkowski
inequality, the space $S_{q'}(S_p)$ should have $\Sig$-cotype
$q'$. However, recalling that $$S_{q'}(S_p) = C_{q'} \otimes_h C_p
\otimes_h R_p \otimes_h R_{q'},$$ we conclude that the subspace
$\mathcal{C}(p,q) = C_{q'} \otimes_h C_p = C_{q'} \otimes_h
R_{p'}$ of $S_{q'}(S_p)$ must also have $\Sig$-cotype $q'$. Then,
we proceed as in Theorem \ref{Theorem-Sigma-Type}. Namely, let us
consider a function $f: \Omega \to \mathcal{C}(p,q)$ so that
$$\widehat{f}_{\mathbf{S}}(\xi) = 0 \qquad \mbox{for} \qquad \xi
\in \Sig \setminus \{\sig\}$$ and such that
$$\widehat{f}_{\mathbf{S}}(\sig) = \sum_{i,j=1}^{d_{\sig}} e_{i1}
\otimes e_{i1} \otimes e_{1j} \otimes e_{1j} \in C_{q'}^{d_{\sig}}
\otimes_h \mathcal{C}(p,q) \otimes_h R_{q'}^{d_{\sig}}.$$ By the
definition of $\Sig$-cotype we have
$$\|\widehat{f}_{\mathbf{S}}(\sig)\|_{S_{q'}^{d_{\sig}}(\mathcal{C}(p,q))}
\le \mathcal{K}_{q'}^2(\mathcal{C}(p,q),\mathbf{S}_{\Sig}) \,
\Big( \int_{\Omega} \Big\| d_{\sig}^{1/q} \mbox{tr}
(\widehat{f}_{\mathbf{S}}(\sig) \zeta^{\sig})
\Big\|_{\mathcal{C}(p,q)}^q d\mu \Big)^{1/q}.$$ Now using the
Banach space isometry $$S_s^n = C_u^n \otimes_h R_v^n \qquad
\mbox{for} \qquad \frac{1}{s} = \frac{1}{2u} + \frac{1}{2v},$$
which follows easily by complex interpolation, we obtain
$$\|\widehat{f}_{\mathbf{S}}(\sig)\|_{S_{q'}^{d_{\sig}}(\mathcal{C}(p,q))}
= \Big\| \sum_{k=1}^n e_{kk} \Big\|_{C_{q'}^{d_{\sig}} \otimes _h
C_{q'}^{d_{\sig}}} \Big\| \sum_{k=1}^n e_{kk}
\Big\|_{C_p^{d_{\sig}} \otimes _h R_{q'}^{d_{\sig}}} =
d_{\sig}^{\frac{1}{2}+ \frac{1}{2p} + \frac{1}{2q'}}.$$ Moreover,
since $\mbox{tr}(\widehat{f}_{\mathbf{S}}(\sig)
\zeta^{\sig}(\omega)) = \sum_{i,j = 1}^{d_{\sig}} e_{ij} \otimes
\zeta_{ji}^{\sig}(\omega) = (\zeta^{\sig}(\omega))^{\mbox{t}}$, we
have $$\Big( \int_{\Omega} \Big\| d_{\sig}^{1/q} \mbox{tr}
(\widehat{f}_{\mathbf{S}}(\sig) \zeta^{\sig}(\omega))
\Big\|_{\mathcal{C}(p,q)}^q d\mu(\omega) \Big)^{1/q} =
d_{\sig}^{\frac{1}{q} + \frac{1}{2q'} + \frac{1}{2p'}}.$$
Combining our previous results, we obtain $q' \ge p'$. This
completes the proof. \fin

\begin{remark}
\emph{By duality, the sharp $\Sig$-type index of $L_p$ is
$\min(p,p')$.}
\end{remark}

\begin{remark}
\emph{By a standard argument using the contraction principle, our
results for sharp $\Sig$-type and $\Sig$-cotype also hold for
\emph{any} uniformly bounded quantized orthonormal system. The
reader is referred to \cite{GP2} for further details.}
\end{remark}

\begin{remark}
\emph{As it was recalled in Remark \ref{Remark-Non-Trivial}, it
seems that there is no analog of the Maurey-Pisier theorem for
operator spaces. Theorem \ref{Theo-Intro-Sharp} clearly reinforces
that idea. Finally, the reader is referred to Section 4.2 of
\cite{J1} for an unrelated notion of operator space cotype $2$ for
which $L_p$ has cotype $2$ whenever $1 \le p \le 2$.}
\end{remark}

\bibliographystyle{amsplain}

\begin{thebibliography}{10}
\bibitem{B} A. Beck, \emph{A convexity condition in Banach
spaces and the strong law of large numbers}, Proc. Amer. Math.
Soc., \textbf{13} (1962), 329-334.
\bibitem {GMP} J. Garc\'{\i}a-Cuerva, J.M. Marco and J. Parcet,
\emph{Sharp Fourier type and cotype with respect to compact
semisimple Lie groups}, Trans. Amer. Math. Soc. \textbf{355}
(2003), 3591-3609.
\bibitem {GP1} J. Garc\'{\i}a-Cuerva and J. Parcet,
\emph{Vector-valued Hausdorff-Young inequality on compact groups},
Proc. London Math. Soc. To appear.
\bibitem {GP2} J. Garc\'{\i}a-Cuerva and J. Parcet, \emph{Quantized
orthonormal systems: A non-commutative Kwapie\'n theorem}, Studia
Math. \textbf{155} (2003), 273-294.
\bibitem {G} D.P. Giesy, \emph{On a convexity condition in
normed linear spaces}, Trans. Amer. Math. Soc. \textbf{125}
(1966), 114-146.
\bibitem {HL} A. Hinrichs and H.H. Lee, \emph{Duality of Fourier
type with respect to locally compact abelian groups}. Preprint.
\bibitem {J1} M. Junge, \emph{Factorization Theory for Spaces of
Operators}. Habilitation Thesis. Kiel, 1996.
\bibitem {J3} M. Junge, \emph{Embeddings of non-commutative
$L_p$-spaces in non-commutative $L_1$-spaces, $1 < p < 2$}, GAFA
\textbf{10} (2000), 389-406.
\bibitem {J2} M. Junge, \emph{Doob's inequality for non-commutative
martingales}, J. reine angew. Math. \textbf{549} (2002), 149-190.
\bibitem {J4} M. Junge, \emph{Embedding of the operator space
$\mathrm{OH}$ and the logarithmic \lq little Grothendieck
inequality\rq}. Preprint.
\bibitem {JX} M. Junge and Q. Xu, \emph{Noncommutative
Burkholder/Rosenthal inequalities}, Ann. Probab. \textbf{31}
(2003), 948-995.
\bibitem {Le} H.H. Lee, \emph{Sharp type and cotype with respect
to quantized orthonormal systems}. Preprint.
\bibitem {Lu} F. Lust-Piquard, \emph{In\'{e}galit\'{e}s de Khintchine dans
$C_p$ $(1 < p < \infty)$}, C.R. Acad. Sci. Paris \textbf{303}
(1986), 289-292.
\bibitem {LuP} F. Lust-Piquard and G. Pisier,
\emph{Non-commutative Khintchine and Paley inequalities}, Ark.
Mat. \textbf{29} (1991), 241-260.
\bibitem {MP} M.B. Marcus and G. Pisier, \emph{Random
Fourier Series with Applications to Harmonic Analysis}, Annals of
Math. Studies \textbf{101}, Princeton University Press, 1981.
\bibitem {MauP} B. Maurey and G. Pisier, \emph{S\'{e}ries de variables
al\'{e}atoires vectorielles ind\'{e}pendantes et propri\'{e}t\'{e}s g\'{e}om\'{e}triques
des espaces de Banach}, Studia Math. \textbf{58} (1976), 45-90.
\bibitem {OR} T. Oikhberg and H.P. Rosenthal, \emph{On certain
extension properties for the space of compact operators}. J. Func.
Anal. To appear.
\bibitem {Pa} J. Parcet, \emph{$\mathrm{B}$-convex operator
spaces}, Proc. Edinburgh Math. Soc. \textbf{46} (2003), 649-668.
\bibitem {Pa2} J. Parcet, \emph{An\'{a}lisis arm\'{o}nico no conmutativo y
geometr\'{\i}a de espacios de operadores}. Ph.D. Thesis. Madrid, 2003.
\bibitem{Pi} G. Pisier, \emph{Sur les espaces qui ne
contiennent pas de $l_1^n$ uniform\'{e}ment}, S\'{e}minaire
Maurey-Schwartz 1973/74, Expos\'{e} 7, \'{E}cole Polytechnique, Paris,
1975.
\bibitem {P0} G. Pisier, \emph{Holomorphic semi-groups and the
geometry of Banach spaces}, Annals of Math. \textbf{115} (1982),
375-392.
\bibitem {P2} G. Pisier, \emph{Non-Commutative Vector Valued
$L_p$-Spaces and Completely $p$-Summing Maps}, Ast\'{e}risque (Soc.
Math. France) \textbf{247} (1998).
\bibitem {P3} G. Pisier, \emph{Introduction to
Operator Space Theory}, Cambridge Univ. Press, 2003.
\bibitem {PS} G. Pisier and D. Shlyakhtenko, \emph{Grothendieck's
theorem for operator spaces}, Invent. Math. \textbf{150} (2002),
185-217.
\bibitem {PX} G. Pisier and Q. Xu, \emph{Non-commutative
martingale inequalities}, Comm. Math. Phys. \textbf{189} (1997),
667-698.
\bibitem {TJ} N. Tomczak-Jaegermann, \emph{The moduli of
smoothness and convexity and the Rademacher ave\-rages of trace
classes $S_p$}, Studia Math. \textbf{50} (1974), 163-182.
\bibitem {X} Q. Xu, \emph{Interpolation of operator spaces}, J. Func.
Anal. \textbf{139} (1996), 500-539.
\end{thebibliography}

\end{document}